\documentclass[11pt,a4paper,psamsfonts]{amsart}
\usepackage{amssymb,amscd,amsxtra,calc}
\usepackage{cmmib57}

\usepackage{graphicx} 
\usepackage{booktabs}
\usepackage{ragged2e}
\usepackage{float}
\usepackage{soul}
\usepackage[citecolor=Navy]{hyperref}
\usepackage[ruled,vlined]{algorithm2e}
\usepackage{amsfonts,amssymb,amsmath,amsthm}
\usepackage{algorithmic}
\usepackage{amssymb,amscd,amsxtra,calc}
\usepackage{cmmib57}
\usepackage{graphicx}
\usepackage{mathrsfs,enumerate,xcolor}

\usepackage[all]{xy}
\usepackage{multirow} 
\usepackage{psfrag,xmpmulti,amscd,color,pstricks, import}
\usepackage{tikz-cd}

\usepackage{graphicx}

\usepackage{lipsum}
\usepackage{amsfonts}
\usepackage{graphicx}
\usepackage{epstopdf}
\usepackage{algorithmic}
\usepackage{float}
\restylefloat{table}
\usepackage{placeins}
\usepackage{array}
\ifpdf
  \DeclareGraphicsExtensions{.eps,.pdf,.png,.jpg}
\else
  \DeclareGraphicsExtensions{.eps}
\fi

\usepackage{subcaption}

\theoremstyle{plain}
    \newtheorem{thm}{Theorem}[section]

    \newtheorem{lemma}[thm]{Lemma}
    \newtheorem{proposition}[thm]{Proposition}

    \newtheorem{theorem}[thm]{Theorem}
\newtheorem{example}[thm]{Example}

\theoremstyle{definition}
    \newtheorem{definition}[thm]{Definition}

\theoremstyle{remark}

\title[Algorithms on Riemannian manifolds and Banach spaces with global guarantees]{Some iterative algorithms on Riemannian manifolds and Banach spaces with good global convergence guarantee}

\author{Tuyen Trung Truong}
\address{Matematikk Institut, Universitetet i Oslo, Blindern, 0851 Oslo, Norway}
\email{tuyentt@math.uio.no}
\thanks{}
\date{\today}


\begin{document}
\begin{abstract}

In this paper, we introduce some new iterative optimisation algorithms on Riemannian manifolds and Hilbert spaces which have good global convergence guarantees to local minima. More precisely, these algorithms have the following properties: If $\{x_n\}$ is a sequence constructed by one such algorithm then: 
 
- Finding critical points:  Any cluster point of $\{x_n\}$ is a critical point of the cost function $f$.

- Convergence guarantee: Under suitable assumptions, the sequence $\{x_n\}$ either converges to a point $x^*$, or diverges to $\infty$. 

- Avoidance of saddle points: If $x_0$ is randomly chosen, then the sequence $\{x_n\}$ cannot converge to a saddle point. 

Our results apply for quite general situations: the cost function $f$ is assumed to be only $C^2$ or $C^3$, and either $f$ has at most countably many critical points (which is a generic situation) or satisfies certain Lojasiewicz gradient inequalities. To illustrate the results, we provide a nice application with optimisation over the unit sphere in a Euclidean space. 

As for tools needed for the results, in the Riemannian manifold case we introduce a notion of "strong local retraction" and (to deal with Newton's method type) a notion of "real analytic-like strong local retraction". In the case of Banach spaces, we introduce a slight generalisation of the notion of "shyness", and design a new variant of Backtracking New Q-Newton's method which is more suitable to the infinite dimensional setting (and in the Euclidean setting is simpler than the current versions).

\end{abstract}
\maketitle

 \section*{Author's information} Department of Mathematics, University of Oslo, Norway. Email: tuyentt@math.uio.no. ORCID: 0000-0001-9103-0923. 

 \section*{Acknowledgements} This paper is based on the author's preprints arXiv: 2001.05768 and arXiv:2008.11091, and supersede them. Here are main new features in this paper compared to these preprints: 
 
 - We focus exclusively on the settings of Riemannian manifolds and Hilbert spaces.  The expositions and proofs in this paper for the Riemannian manifolds are rewritten to help with comprehensibility. 
 
 - In the infinite dimensional setting,  the preprint arXiv: 2001.05768 treated a variant of Backtracking Gradient Descent on Banach spaces. However, the results there implicitly use an assumption that the concerned maps have factors on finite dimensional subspaces, which is very restrictive. In this paper, we treat more general maps over Banach spaces, introducing the following two new objects: 
 
 i) A slight generalisation of the notion of "shyness" \cite{christensen}\cite{hunt-etal}. 
 
ii)  A new variant of Backtracking New Q-Newton's method which is more suitable for the infinite dimensional setting. In the Euclidean setting, this new version is simpler than the current version in  \cite{truong1}. 

 The author was partially supported by Young Research Talents grant 300814 from Research Council of Norway.  We would like to thank Torus AI, as well as Hang-Tuan Nguyen, for help and inspiring discussions.   
 
 \section*{Statements and Declarations} There are no competing interests. 

\section{Introduction} This paper concerns unconstrained optimisation on Riemannian manifolds and Banach spaces. The setting is as follows. Let $X$ be a space (either a Riemannian manifold or a Banach space), and $f:X\rightarrow \mathbf{R}$. We want to solve the problem: $\arg\min _{x\in X}f(x)$. This problem contains, as a special case, the problem of finding roots of a system of equations on $X$: If we want to find a root of a system $F_1(x)=\ldots =F_N(x)=0$, then we consider the function $f(x)=F_1(x)^2+\ldots +F_N(x)^2$, for which roots of the original system are global minima. 

This problem can be helpful to many {\bf constrained} optimisation problems on Euclidean spaces: if the constrained set $S\subset \mathbf{R}^m$ is a Riemannian manifold, then the optimization problem $\arg \min _{x\in S}f(x)$ might be better treated as a Riemannian optimisation problem on the manifold $S$ (where the gradient and Hessian matrix of $f$ is calculated with respect to the induced Riemannian metric on $S$) than the original constrained optimisation problem in the Euclidean space $\mathbf{R}^m$. In the applications later in this paper, we illustrate this with the case where $S$ is the unit sphere. We introduce in this paper some new notations in order to achieve good global convergence guarantee results: "strong local retraction" and (to deal with Newton's method type) a notion of "real analytic-like strong local retraction". 

Another case treated in this paper is that of optimisation problems on Banach spaces. This is motivated by Calculus of variations, where a Partial Differential equation (PDE) can be reduced to minimising a functional $f$ defined on an infinitely dimensional Banach space. The infinite dimensionality presents new difficulties which do not happen in finite dimensions, for example a bounded sequence may not has any convergent subsequence in the strong topology. To this end, we resolve by a reduction to finite dimensions, as well as utilise the notion of "shyness" to measure randomness. 

Since it will need considerable setups  before the new algorithms and results can be discussed unambiguously on Riemannian manifolds and Banach spaces, we defer to Sections 2 and 3 for detail (including some overview of existing work in these settings). In the remaining of this section, we define precisely what we mean to be algorithms with good global convergence guarantees to (local) minima, versions of our algorithms in the more familiar setting of Euclidean spaces, and state the main results.  
 
Given that a closed form formula for solutions to an optimisation problem is in general nonexistent, it is both practically and theoretically the need to make use of an iterative optimisation algorithm. An iterative optimisation algorithm is a predetermined rule $R$ which, starting with a (random) initial point $x_0\in X$, to construct a sequence $x_{n+1}=R(x_n)$. The hope is that the sequence $\{x_n\}$ will converge to a global minimum $x^*$ of $f(x)$. 

Thus, optimisation is a good playground for dynamical systems methods. Indeed, the iterative behaviour of Newton's method for finding roots of a polynomial in 1 complex variable was an origin for the field of Complex Dynamics and that topic is still extensively studied nowadays \cite{RefM}. The tools used to prove the main results in this paper include those from Dynamical Systems and Complex Analysis. Also, randomness plays a crucial role in the proofs. 

We need to restrict our goal though. Because finding a global minimum is in general NP-hard, we at most can only hope to prove that the sequence $\{x_n\}$ will converge to a local minimum $x^*$. A local minimum of $f$ must be a critical point of $f$, i.e. must solve the equation $\nabla f(x)=0$.  We recall that a critical point $x^*$ of $f$ is a saddle point if the Hessian matrix $\nabla ^2f(x^*)$ is invertible and has both positive and negative eigenvalues. It is convenient to also consider a more general type of saddle points, which we will call generalised saddle point: A critical point $x^*$ is a generalised saddle point if $\nabla ^2f(x^*) $ has at least one negative eigenvalue. If $x^*$ is a critical point of $f$ and is not a local minimum, then informally it should be a (generalised) saddle point of $f$.

With the above preparation,   we can divide the question of finding local minima of a function by using an iterative algorithm $R$ into two tasks: 

{\bf Task 1: Critical point finding and Convergence guarantee.}

MainTask 1.1: Show that for every initial point $x_0$, any cluster point of the sequence $\{x_n\}$ constructed by the rule $R$ will be a critical point of the function $f$. That is, if a subsequence $\{x_{n_k}\}$ converges to a point $x^*$, then $\nabla f(x^*)=0$. 

SubTask 1.2: Show that under certain conditions, whenever $\{x_n\}$ has a bounded subsequence, then it converges. 

{\bf Task 2: Avoidance of saddle points and Local minima finding guarantee.}  

MainTask 2.1: Show that if the initial point $x_0$ is randomly chosen and $\{x_n\}$ converges to a point $x^*$, then $x^*$ is not a (generalised) saddle point. 
  
SubTask 2.2: Show that under certain conditions, whenever $\{x_n\}$ converges to a point $x^*$, then $x^*$ is a local minimum.

The reason why Task 2 is important is as follows. It has been shown in  \cite{bray-dean} that generically the ratio between minima and other types of critical points becomes exponentially small when the dimension $k$ increases.  Thus, for an iterative method to be useful in general, it should be able to avoid saddle points. Indeed, when one choose a random initial point $z_0$, it is very likely (in particular, in higher dimensions) that $z_0$ is closer to saddle points than to local minima. If the iterative method has the tendency to converge to critical points closer to the initial point $z_0$, then most likely one will end up with a saddle point and not with a local minimum as desired. 

This paper addresses both tasks for specially designed algorithms on Riemannian manifolds and Banach spaces. We next provide definitions of these algorithms in the familiar setting of Euclidean spaces, so that the readers can grasp their essence, even though the precise definitions need to be waited until Sections 2 and 3 where the necessary backgrounds have been fully developed.  

\subsection{A new variant of Backtracking Gradient Descent method.}\label{SubsectionGD}  We first introduce the well known Gradient Descent method (GD). The general version of this method, invented by Cauchy in 1847 \cite{cauchy}, is as follows. Let $\nabla f(x)$ be the gradient of $f$ at a point $x$, and $||\nabla f(x)||$ its Euclidean norm in $\mathbb{R}^k$.  We choose randomly a point $x_0\in \mathbb{R}^k$ and define a sequence
\begin{eqnarray*}
x_{n+1}=x_n-\delta (x_n) \nabla f(x_n),
\end{eqnarray*}
where $\delta (x_n)>0$ (learning rate), is appropriately chosen. We hope that the sequence $\{x_n\}$ will converge to a (global) minimum point of $f$. 

To boost convergence of Gradient Descent method, a variant is Backtracking Gradient Descent (Backtracking GD) which which works as follows. We fix real numbers $\delta _0>0$ and $0<\alpha ,\beta <1$. We choose $\delta (x_n)$ to be the largest number  $\delta $ among the sequence $\{\beta ^m\delta _0:~m=0,1,2,\ldots\}$ satisfying the Amijo's condition:
\begin{eqnarray*}
f(x_n-\delta \nabla f(x_n))-f(x_n)\leq -\alpha \delta ||\nabla f(x_n)||^2.  
\end{eqnarray*}

It is known, see \cite{truong-nguyen, truong-nguyen2} that Backtracking GD satisfies requirements of Task 1, but it is unknown if it is also so for Task 2 (even though it is seen in many experiments that it seems to be so). This is the reason we introduce the following new version of Backtracking Gradient Descent (Backtracking GD) method: 

\begin{definition}
(Local Backtracking GD.) Let $f:\mathbb{R}^k\rightarrow \mathbb{R}$ be a $C^1$ function. Assume that there are {\bf continuous} functions $r,L:\mathbb{R}^k\rightarrow (0,\infty)$ so that for each $x\in \mathbb{R}^k$, the map $\nabla f$ is Lipschitz continuous on $B(x,r(x))$ with Lipschitz constant $L(x)$. 

The Backtracking GD-Local version procedure is defined as follows. Fix $\delta _0>0$ and $0<\alpha , \beta <1$.  For each $x\in \mathbb{R}^k$, we define $\widehat{\delta} (x)$ to be the largest number $\delta$ among  $\{\beta ^n\delta _0:~n=0,1,2,\ldots \}$ which satisfies the two conditions
\begin{eqnarray*}
\delta &<&\alpha /L(x),\\
\delta ||\nabla f(x)|| &<& r(x).  
\end{eqnarray*}

For any $x_0\in \mathbb{R}^k$, we then define the sequence $\{x_n\}$ as follows
\begin{eqnarray*}
x_{n+1}=x_n-\widehat{\delta} (x_n)\nabla f(x_n). 
\end{eqnarray*}
\label{DefinitionBacktrackingGDNew}\end{definition}

For later use, here we recall that a function $f$ is in class $C^{1,1}_L$, for a positive number $L<\infty$, if $\nabla f(x)$ is {\bf globally} Lipschitz continuous with Lipschitz constant $L$. We also recall that a point $x_{\infty}$ is a cluster point of a sequence $\{x_n\}$ if there is a subsequence $\{x_{n_j}\}$ which converges to $x_{\infty}$.

{\bf Examples.}  (i) If $f$ is in $C^{1,1}_L$, then $f$ satisfies the condition in Definition \ref{DefinitionBacktrackingGDNew} by defining $L(x)=L$ for all $x$.  (ii) If $f$ is in $C^2$, then we can choose any continuous function $L(x)$ so that $L(x)\geq \max _{z\in B(x,r(x))}||\nabla ^2f||$, where $r:\mathbb{R}^k\rightarrow (0,\infty )$ is any continuous function.

{\bf Motivation.} Here is a motivation for the Local Backtracking GD algorithm. If $\nabla f$ is Lipschitz continuous, then it is known (see \cite{armijo}) that the largest possible choice of the learning rate so that Backtracking GD has global convergence guarantee is $\alpha /L$. For a general cost function, the Local Backtracking GD uses a local Lipschitz continuity estimate for $\nabla f(x)$. Theorem 3.5 in \cite{truong-nguyen} shows that roughly speaking (and precisely, when the dimension $k=1$), the choice of the learning rate in Local Backtracking GD is the largest one can in order to always guarantee global convergence. 

{\bf Practical implementation.} For ease of implementation, we can choose $r(x)$ to be a constant, for example $\epsilon _0$ for a small positive number $\epsilon _0$, and estimate $L(x)$ as $||\nabla ^2f(x)||+\gamma _0 (x,\epsilon _0)$ by a simple Taylor's expansion (for example, $\gamma _0(x,\epsilon )=2\epsilon _0||\nabla ^3f(x_0)||$). Alternatively, we can simply choose a small positive constant $\gamma _0>0$ and choose the learning rate as $\delta =\alpha /(||\nabla ^2f(x)||+\gamma _0)$.  

{\bf Remark.} In the infinite dimensional setting, a similar algorithm can be used. 

\subsection{(Backtracking) New Q-Newton's method} Among higher order optimisation methods, the most well known one is probably Newton's method. Assume that $f$ is a $C^2$ function. Then Newton's method constructs a sequence $\{z_n\}$ by the following update rule: $$z_{n+1}=z_n-(\nabla ^2f(z_n))^{-1}.\nabla f(z_n).$$ 

An attractive feature of Newton's method is its quick rate of convergence (when it actually converges). Also, it is straight forward to implement.  On the other hand, it has many disadvantages. First, it only works if $\nabla ^2f(z_n)$ is invertible. Second, it has the tendency to converge to the nearby critical point even if that point is a saddle point (this can be readily checked with a quadratic function), and hence fails Task 2. Third, its implementation in large scale can be costly. In this paper, we concentrate only the first two issues.

In this paper, we study more deeply the algorithm New Q-Newton's method from \cite{truong-etal}. 

Let $A:\mathbb{R}^k\rightarrow \mathbb{R}^k$ be an invertible {\bf symmetric} square matrix. In particular, it is diagonalisable.  Let $V^{+}$ be the vector space generated by eigenvectors of positive eigenvalues of $A$, and $V^{-}$ the vector space generated by eigenvectors of negative eigenvalues of $A$. Then $pr_{A,+}$ is the orthogonal projection from $\mathbb{R}^k$ to $V^+$, and  $pr_{A,-}$ is the orthogonal projection from $\mathbb{R}^k$ to $V^-$. As usual, $Id$ means the $k\times k$ identity matrix.

\medskip
{\color{blue}
 \begin{algorithm}[H]
\SetAlgoLined
\KwResult{Find a minimum of $F:\mathbb{R}^k\rightarrow \mathbb{R}$}
Given: $\{\delta_0,\delta_1,\ldots, \delta_{k}\}\subset \mathbb{R}$\  and $\alpha >0$;\\
Initialization: $x_0\in \mathbb{R}^k$\;
 \For{$n=0,1,2\ldots$}{ 
    $j=0$\\
    \If{$\|\nabla f(x_n)\|\neq 0$}{
   \While{$\det(\nabla^2f(x_n)+\delta_j \|\nabla f(x_n)\|^{1+\alpha}Id)=0$}{$j=j+1$}}

$A_n:=\nabla^2f(x_n)+\delta_j \|\nabla f(x_n)\|^{1+\alpha}Id$\\
$v_n:=A_n^{-1}\nabla f(x_n)=pr_{A_n,+}(v_n)+pr_{A_n,-}(v_n)$\\
$w_n:=pr_{A_n,+}(v_n)-pr_{A_n,-}(v_n)$\\
When $f$ does not have compact sublevels, normalise $w_n:=w_n/\max\{1,||w_n||\}$\\
$x_{n+1}:=x_n-w_n$
   }
  \caption{New Q-Newton's method} \label{table:alg}
\end{algorithm}
}
\medskip

To boost convergence, \cite{truong1} incorporated two new ideas: 

1. Choose $\delta _j$ so that all eigenvalues of the matrix $\nabla^2f(x_n)+\delta_j \|\nabla f(x_n)\|^{1+\alpha}Id$ has absolute value $\geq $ $\|\nabla f(x_n)\|^{1+\alpha} \min _{i\not= j}|\delta _i-\delta _j|/2$ (which is possible by the Pigeon hole principle). 

2. Armijo's Backtracking line search in to New Q-Newton's method (which is possible by noting that the inner product $<w_n,\nabla f(x_n)>$ is non-negative), to obtain a new algorithm named Backtracking New Q-Newton's method:
$$w_{n+1}=x_n-\gamma _nw_n,$$  
where $\gamma _n$ is the learning rate found by using Armijo's Backtracking line search. 

{\bf Remark.} In the infinite dimensional setting, the algorithm above poses some huge challenges. We therefore later will design a new version more suitable for that setting. In the Euclidean setting, that version turns out to be simpler than the current version. 

\subsection{The main results}\label{SubsectionMainResults} We are now ready to state the main results of this paper.  We recall that the Kurdyka-Lojasiewicz gradient inequality  condition is satisfied at a point $z^*$ if there is an open neighbourhood $U$ of $z^*$ and constants $C>0$, $0<\mu <1$ such that for all $z\in U$ we have $|f(z)-f(z^*)|^{\mu}\leq C||\nabla f(z)||$.  This is an important class of cost functions to which there is intensive research, in particular in complex variables, real algebraic geometry and algebraic optimisation.  

\begin{theorem} Let $X$ be a Riemannian manifold, and $f:X\rightarrow \mathbf{R}$ an objective function of class $C^1$ (in case we are using the first order methods) or $C^2$ (in case we are using second order methods).

1) For Critical point finding and Convergence guarantee: 

1i) Both Backtracking GD, Local Backtracking GD, and Backtracking New Q-Newton's method satisfy the requirements of MainTask 1.1.

1ii) If $f$ has at most countably many critical points, then both  Backtracking GD, Local Backtracking GD, and Backtracking New Q-Newton's method satisfy the requirements of SubTask 1.2. 

1iii) If $f$ satisfies the Kurdyka-Lojasiewicz gradient inequality condition at its critical points, then both Backtracking GD and Local Backtracking GD satisfy the requirements of SubTask 1.2. 

2) For Avoidance of saddle points and Local minima finding guarantee: 

2i) If $f$ is $C^2$ near its generalised saddle points and the parameters are randomly chosen, then Local Backtracking GD satisfies requirements of MainTask 2.1. 

2ii) If $f$ is $C^3$ near its saddle points and the parameters are randomly chosen, then both New Q-Newton's method and Backtracking New Q-Newton's method satisfy requirements of MainTask 2.1.  

\label{TheoremRiemannianManifold}\end{theorem}

This result will be proven in Section 2 (separately for first and second order methods), where precise conditions on the randomness of parameters (which are complicated to state here) are given in detail. Note that Local Backtracking GD also satisfies a somewhat stronger requirement than MainTask 2.1:  we can show that {\bf any cluster point} of the constructed sequence $\{x_n\}$ cannot be an {\bf isolated} generalised saddle point.  

A similar result to Theorem \ref{TheoremRiemannianManifold} is proven also for the case where $X$ is a Banach space, see Section 3, where necessary modifications (e.g. weak convergence instead of strong convergence) are needed. To overcome the difficulty of having infinite dimensions, we introduce some new ideas: defining a slight generalisation of the "shyness" \cite{christensen}\cite{hunt-etal} and designing a new variant of Backtracking New Q-Newton's method. More precisely, we show that Local Backtracking GD solves both Tasks 1 and 2 from the introduction, while BNQN solves Task 1 and locally solves Task 2 (i.e. local Stable-center manifold exists for the dynamics of BNQN near saddle points, and BNQN has quadratic convergence rate near non-degenerate local minima). To be able to show global avoidance of saddle points for BNQN in the Banach space setting, it is necessary to extend Theorem \ref{TheoremRandomnessLambda} below to the infinite dimensional setting, which is beyond the current paper. 

\section{Optimisation/Root-finding on Riemannian manifolds} Many constrained optimisation on Euclidean spaces have the eligible set to be a Riemannian manifold. As such we can consider them in the most general setting of optimisation on Riemannian manifolds. We note that this setting also offers many interesting and useful applications, such as in the Netflix prize competition \cite{netflix}. Another example is Fisher information metric, which is a natural choice of Riemannian metric on a space of parametrised statistical models (such as Deep Neural Networks), the modern theory is largely due to S. Amari \cite{amari-nagaoka}. Yet one other interesting case is that of constrained optimisation problems on Riemannian manifolds, if the constraints give rise naturally to a Riemannian submanifold. 

Since the algorithms involve several important notions from Riemannian geometry, we first present some preliminaries on them. After that, we will give precise definitions of the new algorithms, and then proceed with main results, their proofs and applications. 

\subsection{Preliminaries on Riemannian geometry}\label{SubsectionPreliminaries} This is just a very terse overview of backgrounds from Riemannian geometry needed for later use. We refer the interested readers to \cite{lee} (for generalities about Riemannian manifolds) and \cite{absil-etal, boumal} (for more details on how to use the tools for optimisation on Riemannian manifolds). 

A Riemannian manifold is a manifold $X$, together with a Riemannian metric $g(x)$ which is an inner product on tangent spaces $T_xX$. Usually we assume that $g(x)$ varies smoothly when $x$ changes. With the help of a Riemannian metric, and the associated Levi-Civita connection, given a function $f:X\rightarrow \mathbb{R}$ we can define the notions of gradient $grad(f)\in TX$ and Hessian $Hess(f)$ which reduce to the familiar notions $\nabla f$ and $\nabla ^2f$ when $X$ is a Euclidean space. 

A Riemannian metric gives rise to a metric $d_X$ on $X$ in the following manner. If $\gamma :[a,b]\rightarrow X$ is a smooth curve, then we define the length of $\gamma$ as:

\begin{eqnarray*}
L(\gamma )=\int _a^b||\gamma '(t)||dt.
\end{eqnarray*}
Here $||\gamma '(t)||$ is the length of the vector $\gamma '(t)$, with respect to the inner product given by the Riemannian metric $g$. 

Given $x,y\in X$, the distance $d_X(x,y)$ is the infimum  of the  lengths $L(\gamma )$, where $\gamma $ runs on all over curves on $X$ connecting $x$ and $y$. 

A geodesic is a curve which realises the distance. It satisfies a second order ODE. By results from ODE, uniqueness and local existence of geodesics are guaranteed. In particular, one can define exponential maps. If $x\in X$ and $v\in T_xX$ with $||v||$ small enough, then there is a unique geodesic $\gamma _v:[0,1]\rightarrow X$ such that $\gamma _v(0)=x$ and $\gamma _v '(0)=v$. The exponential map is $exp _x(v)=\gamma _v(1)$.  

For $x\in X$ and $r>0$, we denote by $B(T_xX,r)$ the set $\{v\in T_x:~||x||<r\}$. The injectivity radius, $inj(x)$, is the supremum of all $r$, for which the exponential map $exp_x$ is well-defined and a diffeomorphism from $B(T_x,r)$ onto its image. By the previous paragraph, we always have $inj (x)>0$. The injectivity radius can be infinity, for example in the case $X$ is a Euclidean space. We note the following important property of injectivity radius (see Proposition 10.18 in \cite{boumal}). 

\begin{proposition} The map $inj:~X\rightarrow (0,\infty ]$ is continuous. 
\label{PropositionInjectivityRadius}\end{proposition}

For a complete Riemannian manifold, a proof can be found in \cite[Proposition 10.37]{lee}. For a general Riemannian manifold, a proof is given by M. Stephen and J. Lee in an online discussion \cite{inj} and is incorporated into  \cite[Section 10.8]{boumal}. 

Exponential maps provide a way to move, on the same manifold, from a point $x$ in a chosen direction $v$. However, in optimisation, the most crucial property of an exponential map is that its derivative at $0\in T_xX$ is the identity map. (This implies, by inverse function theorem, that the exponential map is a local diffeomorphism near $0$.) This is generalised to the following (global) notion of "retraction", as given in \cite[Definition 4.1.1]{absil-etal} which we now recall. 

{\bf Definition (Global retraction).} A retraction on a (Riemannian) manifold $X$ is  a smooth mapping $R:TX\rightarrow X$ from the tangent bundle $TX$, with the following properties. If $R_x=R|_{T_xX}$ then: i) $R_x(0_x)=x$ where $0_x$ is the zero element of $T_xX$, and ii) $DR_x(0_x)=Id_{T_xX}$.  

On Euclidean spaces or complete Riemannian manifolds, global retractions exist. Some other interesting global retractions are given in \cite[Section 4.1]{absil-etal}. On the other hand, it is not clear if all Riemannian manifolds have at least one retraction in the above sense (we think that probably the answer is No). For example, by Hopf-Rinow theorem, the exponential map is defined on the whole tangent bundle if and only if $X$ is a complete metric space. In the current literature, all theoretical results are stated and proven under the existence of such global retractions. There are also a notion of "local retractions" \cite[Section 4.1.3 ]{absil-etal}, but we are not aware of any use of them in theoretical treatments, the reasons may be that they are not strong enough. In the next subsection, we will discuss a stronger version called  "strong local retractions", which exist on all Riemannian manifolds and which are strong enough to guarantee good theoretical properties.    

Next we discuss some estimates for Taylor's expansion of functions of the form $f(R_x(v))$ where $R_x:~B(T_xX,r(x))\rightarrow X$ a diffeomorphism onto its image such that $R_x(0)=x$ and $DR_x(0)=Id_{T_xX}$. The presentation here is taken from \cite[Section 10.4]{boumal}. We let $\widehat{f_x}=f\circ R_x:B(T_xX,r(x))\rightarrow \mathbb{R}$.  Suppose that 
\begin{equation}
||\nabla \widehat{f_x}(v)-\nabla \widehat{f_x}(w) ||\leq L||v-w||,
\label{Equation1}\end{equation}
for all $v,w\in B(T_xX,r)$, here $L>0$ is a positive constant. Then, 
\begin{equation}
|f(R_x(v))-f(x)-<v,grad f(x)>|\leq L||v||^2/2,
\label{Equation2}\end{equation}
 for all $v\in B(T_xX,r)$. This inequality is the generalisation of the usual property of functions on an open subset of a Euclidean space whose gradient is (locally) Lipschitz continuous. While Equation (\ref{Equation1})  can be complicated for general maps $R_x$ and functions $f$, there is one special case where it has the usual form. Indeed, if $f$ is $C^2$ and $R_x=exp_x$ is the exponential map and $r(x)=inj(x)$, then Equation (\ref{Equation1}) is satisfied for all $x\in X$ if and only if $||Hess (f)||\leq L$ on $X$.

A more general inequality is as follows, see \cite[Exercise 10.51]{boumal}. Let notations be as in the previous paragraph. Suppose that 
\begin{equation}
||\nabla \widehat{f_x}(v)-\nabla \widehat{f_x}(0) ||\leq L||v||,
\label{Equation1}\end{equation}
for all $v\in B(T_xX,r)$, here $L>0$ is a positive constant. Then, 
\begin{equation}
|f(R_x(v))-f(x)-<v,grad f(x)>|\leq L||v||^2/2,
\label{Equation2}\end{equation}
 for all $v\in B(T_xX,r)$. 
 
 The above inequalities can be extended globally by the following trick. If $\gamma :[0,1]\rightarrow X$ is a continuous curve, with $x=\gamma (0)$ and $y=\gamma (1)$, then we can find numbers $t_0=0<t_1<\ldots <t_{N-1}<t_N=1$, so that for all $j$, the point $x_j=\gamma (t_j)$ is in the range of the exponential map $exp_{x_{j-1}}$. Moreover $x_{j}=exp_{x_{j-1}}(v_{j-1})$ where $\sum ||v_j||\sim $ the length of $\gamma$. Hence, with the above local inequalities, one can prove global inequalities like $|f(x)-f(y)|\leq Ld_X(x,y)$, by using a telescope sum $\sum |f(x_j)-f(x_{j-1})|$. 

A diffeomorphism $\iota :X\rightarrow \iota (X)\subset \mathbb{R}^N$ is an isometric embedding if the lengths of vectors are preserved. This leads to preservation of many properties, including distances. We end this with Nash's embedding theorem \cite{kuiper, nash1, nash2}. 

\begin{theorem} Let $X$ be a Riemannian manifold of finite dimension. Then there is an isometric embedding $\iota: X\rightarrow \mathbb{R}^N $ for some big enough dimension $N$. 
\label{TheoremNashEmbedding}\end{theorem}

\subsection{(Local-) Backtracking GD and New Q-Newton's method on Riemannian manifolds}\label{SubsectionAlgorithms} The bulk of this subsection is to extend the (Local-) Backtracking GD and New Q-Newton's method  to a general Riemannian manifold, in such a way that known results for optimisation on Euclidean spaces can be extended. 

In the previous subsection, we mentioned that it is not known if global retractions exist on every Riemannian manifolds. On the other hand, as mentioned, the versions of local retractions in the current literature are not strong enough to guarantee good theoretical properties of iterative optimisations. We will first introduce a version of local retractions, called strong local retraction, which is both existing on all Riemannian manifolds and strong enough to derive good theoretical guarantees. We recall that if $x\in X$ and $r>0$, then $B(T_xX,r)$ $=$ $\{v\in T_xX:~||v||<r\}$. 

\begin{definition}

{(Strong Local Retraction).}  A strong local retraction consists of a function $r:X\rightarrow (0,\infty ]$ and a map $R:\bigcup _{x\in X}B(T_xX,r(x))\rightarrow X$ with the following properties: 

1) $r$ is upper semicontinuous, that is $\limsup _{y\rightarrow x}r(y)\leq r(x)$.

2) $\bigcup _{x\in X}B(T_xX,r(x))$ is an open subset of $TX$, and $R$ is $C^1$ on its domain. 

3) If $R_x=R|_{B(T_x,r(x))}$, then $R_x$ is a diffeomorphism and $DR_x(0)=Id$. Moreover, we assume that $R$ is $C^1$.  
\label{Def1}\end{definition}

\begin{example}
1) Because of the assumption that $\bigcup _{x\in X}B(T_xX,r(x))$ is an open subset of $TX$, it follows that for every compact set $K\subset X$ we have $\inf _{x\in K}r(x)>0$. 

2) If $R$ is a global retraction, then with the choice of $r(x)=\infty$ for all $x$, we have that $(r,R)$ is a Strong local retraction. 

3) If $r(x)=inj(x)$ and $R=$ the exponential map, then $(r,R)$ is a Strong local retraction, by Proposition \ref{PropositionInjectivityRadius}. Thus on every Riemannian manifold, there exists a Strong local retraction, where moreover we can assume that $r$ is a continuous function.  
\label{Example1}\end{example}

We have the following useful property on local Lipschitz continuity for Strong local retractions. 

\begin{lemma}
Let $r,R$ be a Strong local retraction on a Riemannian manifold $X$.   For each $x\in X$, there exists $s(x),h(x)>0$ such that the following property holds. For all $y,z\in R(B(T_xX,s(x)))$ and $v\in T_yX$, $w\in T_zX$ so that $||v||<s(x)-||R_x^{-1}(y)||$ and $||w||<s(x)-||R_x^{-1}(z)||$, then 
\begin{eqnarray*}
d_X(R_y(v),R_z(w))\geq \frac{1}{2}d_X(y,z)-h(x)||PT_{x\leftarrow y}v-PT_{x\leftarrow z}w||.
\end{eqnarray*}
Here $d_X$ is the induced metric on $X$, and $PT_{x\leftarrow y}$ is the parallel transport of vectors from $T_yX$ to $T_xX$ along the unique geodesic from $x$ to $y$ (when $y$ is close enough to $x$). 
\label{LemmaLowerBoundR}\end{lemma}
\begin{proof}
Indeed, this is a consequence of the assumption that $R$ is $C^1$ (and hence in particular is locally Lipschitz continuous, in both variables $x$ and $v\in B(T_xX,r(x))$) and the assumption that $R_y(0)=y$ for all $y$, which the readers can ready work out on local coordinate charts. 
\end{proof}

Now we define the update rule for the generalisation of the basic version of Backtracking GD on a general Riemannian manifold.  

\begin{definition}

{(Riemannian Backtracking GD).} Let $X$ be a Riemannian manifold, and $(r,R)$ a Strong local retraction on $X$.  We choose $0<\alpha , \beta <1$ and $\delta _0>0$. Let $f:X\rightarrow \mathbb{R}$ be a $C^1$ function. For $x\in X$, we choose $\delta (x)$ to be the largest number $\delta $ in the set $\{\beta ^j\delta _0:~j=0,1,2,\ldots \}$ which satisfies both $\delta ||grad (f)(x)||< r(x)/2$ and Armijo's condition: 
\begin{eqnarray*}
f(R_x(-\delta grad (f)(x)))-f(x)\leq -\alpha \delta ||grad (f)(x)||^2. 
\end{eqnarray*}
(Since $R_x(0)=0$ and $DR_x(0)=Id$, there exists at least a positive number $\delta '$ for which Armijo's condition is satisfied. Hence, the function $\delta (x)$ is well-defined.)

The update rule for Riemannian Backtracking GD is as follows: We choose an initial point $x_0$, and construct a sequence $x_{n+1}=R_{x_n}(-\delta (x_n)grad (f)(x_n))$. 
\label{Def2}\end{definition}

\begin{example}
In the Euclidean space, this definition is classical, goes back at least to \cite{armijo}. In the Riemannian manifold setting, where the retraction $R$ is global, it is is known for awhile \cite{absil-etal, boumal}. In our setting of Strong local retractions here, it is a bit more complicated to state. 
\label{Example0}\end{example}

\begin{lemma}
Let the setting be as in the definition for the Riemannian Backtracking GD algorithm. Let $K\subset X$ be a compact set. If $\inf _{x\in K}||grad (f)(x)||>0$, then $\inf _{x\in K}\delta (x)>0$. 
\label{LemmaLowerBound}\end{lemma}
\begin{proof}
The proof is exactly as in the Euclidean setting, by using observation 1) in Example \ref{Example1}, which one can find for example in \cite{truong-nguyen}. 
\end{proof}

While it is expected that the basic version of Riemannian Backtracking GD should avoid saddle points (see more detailed discussion at the end of this paper), currently we do not know how to do so. Hence, we next define the update rule for the Riemannian version of Local  Backtracking GD, for which avoidance of saddle points is guaranteed. 

\begin{definition}
{ (Riemannian Local  Backtracking GD)} Let $X$ be a Riemannian manifold, and $(r,R)$ a Strong local retraction on $X$, where it is assumed that $r$ is continuous.  Let $f:X\rightarrow \mathbb{R}$ be $C^1$. Put $\widehat{f_x}=f\circ R_x: B(T_xX,r(x))\rightarrow X$. We assume that there is a continuous function $L:X\rightarrow (0,\infty)$ such that for all $x\in X$ and all $v\in B(T_xX,r(x))$, the following inequality is satisfied: 
\begin{eqnarray*}
||\nabla \widehat{f_x}(v)-\nabla \widehat{f_x}(0)||\leq L(x)||v||. 
\end{eqnarray*}

Moreover, assume that the conclusions of Lemma \ref{LemmaLowerBoundR} are satisfied with the choice $s(x)=r(x)$ and $h(x)=L(x)$.
 
Fix $0<\alpha ,\beta <1$ and $\delta _0$. For each $x\in X$, we define $\widehat{\delta}(x)$ to be the largest number $\delta $ among $\{\beta ^j\delta _0:~j=0,1,2,\ldots \}$ which satisfies the two conditions: 

\begin{eqnarray*}
\delta &<&\alpha /L(x),\\
\delta ||grad (f)(x)|| &<& r(x)/2.
\end{eqnarray*}

The update of Riemannian Local Backtracking GD is as follows. We choose an initial point $x_0\in X$, and construct the sequence $\{x_n\}$ as follows: 
\begin{eqnarray*}
x_{n+1}=R_{x_n}(-\widehat{\delta}(x_n)grad (f)(x_n)). 
\end{eqnarray*}
\label{Def3}\end{definition}

\begin{example} This definition was given in  the introduction in the Euclidean spaces setting, where it is simpler to state.

(i) If $f$ is in $C^{1,1}_L$ (see \cite{absil-etal, boumal} for precise definition in the Riemannian setting, see also the previous subsection), then we can choose $L(x)=L$ for all $x$. 

(ii) If $f$ is $C^2$ and $R=$ the exponential map, then after  making $r(x)$ to be finite (for example, by replacing it with $\min \{r(x),1\}$), we can choose $L(x)$ to be any continuous function so that $L(x)\geq \sup _{z\in R_x(B(T_xX,r(x)))}||Hess(f)||$ for all $x\in X$. 

iii) More generally, if $f,R$ are $C^2$ functions, then since $\widehat{f_x}=f\circ R_x$ is $C^2$, we see that the conditions to apply Riemmanian Local Backtracking GD are fully satisfied. 
\label{Example2}\end{example}

Finally, we define the update rule for the Riemannian version of New Q-Newton's method. 

\begin{definition}
{(Riemannian New Q-Newton's method)} Let $X$ be a Riemannian manifold of dimension $k$ with a Strong local retraction $r,R$. Let $f:X\rightarrow \mathbb{R}$ be a $C^2$ function. We choose a real number $1<\alpha $ and randomly $k+1$ real numbers $\delta _0,\ldots , \delta _k$.  We also choose a strictly increasing sequence $\{\gamma _j\}_{j=0,1,2,\ldots } $, for which $\gamma _0=0$, $\gamma _1=1$, $\lim _{j\rightarrow\infty}\gamma _j=\infty$ and $\liminf _{j\rightarrow\infty}\gamma _j/\gamma _{j+1}>0$.  The update rule for Riemannian New Q-Newton's method is as follows. We choose an initial point $x_0$, and construct the sequence $\{x_n\}$ as follows: 

- If $grad(f)(x_n)=0$, then STOP. Otherwise,

- Choose $j$ to be the smallest number among $\{0,1,\ldots ,k\}$ so that $A_n=Hess(f)(x_n)+\delta _j||grad (f)(x_n)||^{\alpha}Id$ is invertible. 

- Let $V_{A_n}^+\subset T_{x_n}X$ be the vector space generated by eigenvectors with {\bf positive} eigenvalues of $A_n$, and $V_{A_n}^-\subset T_{x_n}X$ be the vector space generated by eigenvectors with {\bf negative} eigenvalues of $A_n$. Let $pr_{\pm,A_n}:T_{x_n}X\rightarrow V_{A_n}^{\pm}$ be the corresponding orthogonal projections.  

- Define $v_n$ by the formula $v_n=A_n^{-1}.grad(f)(x_n)$. 

- Let $w_n=pr_{+,A_n}.v_n-pr_{-,A_n}.v_n$. 

- Choose $j$ to be the smallest number so that $\gamma _j r(x_n)/2\leq v_n<\gamma _{j+1}r(x_n)/2$, then define  $\lambda _n=1/\gamma _{j+1}$, and $x_{n+1}=R_{x_n}(-\lambda _nw_n)$. (If $r(x)=\infty$, then we simply choose $\lambda _n=1$.)

\label{Def4}\end{definition}

As in the Euclidean space setting, here we can also incorporate Armijo's Backtracking line search to obtain the Riemannian Backtracking New Q-Newton's method algorithm. 

\begin{example} This definition was given in  \cite{truong-etal} in the Euclidean setting. There, $\lambda _n=1$ because $r(x_n)=\infty$ for all $x_n$.   

In our definition here, if $||v_n||$ is small (relative to $r(x_n)$), then $\lambda _n=1$. 

The assumption that $\delta _0,\delta _1,\ldots, \delta _k$ are randomly chosen is only needed to make sure that the preimages of sets with {\bf zero} Lebesgue measure, by the dynamical systems associated to Riemannian New Q-Newton's method, are again of {\bf zero} Lebesgue measure (in the literature, this property is sometimes called Lusin ($N^{-1}$) property). This assumption was overlooked in \cite{truong-etal}. On the other hand, experiments in that paper show that the algorithm works well even if we do not choose  $\delta _0,\delta _1,\ldots, \delta _k$ randomly. (See Theorem \ref{TheoremRandomnessLambda} for more detail on the level of randomness required.)  It is possible that indeed Lusin ($N^{-1}$) property can hold for the dynamical systems in Riemannian New Q-Newton's method under much more general assumptions on $\delta _0,\ldots ,\delta _k$.         
\label{Example3}\end{example} 

For a general Riemannian manifold, Riemannian New Q-Newton's method is guaranteed to avoid saddle points {\bf locally}, that is we can establish local Stable-Center manifold theorems near saddle points. Whether this avoidance of saddle points can be established globally for all Riemannian manifolds is to be studied in a future work. In the current paper, we introduce a stronger assumption on the Strong local retraction $r,R$, under which condition global avoidance of saddle points for Riemannian New Q-Newton's method is guaranteed. 

\begin{definition}
{ (Real analytic-like Strong local retraction)} We assume that the Strong local retraction $(r,R)$ on $X$ has the following property. For every point $x_0\in X$, there is a small open neighbourhood $U$ of $x_0$ and a small open set $W\subset \mathbb{R}$ , so that if $u(x,s):~U\times W\rightarrow \bigcup _{x\in U}B(x,r(x))$ is continuous, as well as a $C^1$ map in the variable $x$ and a real rational function in the variable $s$, then for all $y\in U$, there is a real analytic function $h_y$  and a real rational function $\kappa _y$ such that $h_y\circ \kappa _y(s)$ restricts to $\det (grad _y(R_y(u(y,s))))$ on $W$. Here $\det (.)$ is the determinant of a matrix.
\label{Def5}\end{definition}

\begin{example} The use of this condition lies in that if $s\mapsto \det (grad _y(R_y(u(y,s))))$ ($s\in W$)  is not the zero function identically, then its zero set has Lebesgue measure $0$. 

One prototype example for this is when $X$ is an {\bf open subset} of a Euclidean space and $R$ is the exponential map (which is just the identity map in this case). In this case, we just need to take $h\circ \kappa (s)$ to be the rational function which defines the function $s\mapsto  \det (grad _y(R_y(u(y,s))))$ for $s\in W$. 

More generally, this is the case if the Strong local retraction $R(x,v)$ ($x\in X,v\in B(x,r(x))$) can be extended real analytically (need not be globally a diffeomorphism onto its image) to the whole $T_xX$. For example, this is the case when $X$ is an open subset of a  Lie group, and $R_x$ is the exponential map.   
\label{Example4}\end{example}

\subsection{Some useful results about functions defined on Riemannian manifolds}\label{SubsectionMainResults} In this subsection we state and prove main results on convergence and/or avoidance of saddle points results for the algorithms defined in the previous subsection. In addition, we will also prove similar results for a continuous version of Backtracking GD. We will start with some preparation results.  

The following is used to prove properties of the Riemannian Backtracking GD, Riemannian Continuous Backtracking GD as well as Riemannian Local Backtracking GD. We say that a sequence $\{x_n\}$, in a metric space $X$, diverges to infinity if it eventually leaves every compact subsets of $X$. We say that a point $x\in X$ is a cluster point of $\{x_n\}$ if there exists a subsequence $\{x_{n_j}\}$ converging to $x$.  

\begin{theorem} Let $X$ be a Riemannian manifold of finite dimension. Let $d_X$ be the induced metric on $X$. Let $\{x_n\}$ be a sequence in $X$ such that $\lim _{n\rightarrow\infty}d_X(x_{n+1},x_n)=0$. Let $C$ be the set of cluster points of $\{x_n\}$. Let $A$ be a closed subset of $X$. Assume that $C\subset A$. 
Let $B$ be a connected component of $A$, and assume further that $B$ is compact. 

1) Assume that $C\cap B\not=\emptyset$. Then $C\subset B$ and $C$ is connected. 

2) (Capture theorem) Assume that $B$ is a point and $C\cap B\not= \emptyset$. Then $C=B$, i.e. the sequence $\{x_n\}$ converges to the point $B$. 

3) Assume that $A$ is at most countable. Then either $\{x_n\}$ converges to a point, or $\{x_n\}$ diverges to infinity.  

\label{TheoremConvergenceSequenceRiemannianManifolds}\end{theorem}
\begin{proof}
The idea is to apply Nash's embedding theorem to the arguments in \cite{truong-nguyen2, truong-nguyen}.

By Theorem \ref{TheoremNashEmbedding}, for the purpose of comparing metrics, we can assume that $X$ is a Riemannian submanifold of some Euclidean space $\mathbb{R}^N$. We let $||.||$ denote the usual norm on $\mathbb{R}^N$. We let $\mathbb{P}^N$ be the real projective space of dimension $N$, and $d(.,.)$ the standard spherical metric.

Then for $x,x'\in X$, we have the following inequalities: $d_X(x,x')\geq ||x-x'||\geq d(x,x')$. The first inequality follows since $X$ is a Riemannian submanifold of $\mathbb{R}^N$.  The second inequality is probably well known, and a detailed proof is given in \cite{truong-nguyen, truong-nguyen2}.  Hence, we also have $\lim _{n\rightarrow\infty}d(x_{n+1},x_n)=0$. Note that while the metric $d_X$ may be very different from that of the restriction of $d(.,.)$ to $X$ (in particular, since $d(.,.)$ is bounded, while $d_X$ may not), $X$ is a {\bf topological} subspace of $\mathbb{P}^N$. In particular,  convergence behaviour of a sequence $\{x_n\}\subset X$, considered in either the original topology on $X$ or the induced one from $\mathbb{P}^N$, is the same. 

We let $C'\subset \mathbb{P}^N$ be the set of cluster points of $\{x_n\}$, considered as a sequence in $\mathbb{P}^N$. Then $C'$ is the closure of $C$ in $\mathbb{P}^N$, and $C'\cap \mathbb{R}^N=C\subset A$. Since $\mathbb{P}^N$ is a compact metric space, it follows from \cite{asic-adamovic} that $C'$ is a connected set.  

 1) If $C'\cap B=C\cap B$ is not empty, then since $C'\cap \mathbb{R}^N\subset A$, $C'$ is connected and $B$ is a compact connected component of $A$, it follows that $C'\subset B$. Hence, $C'=C'\cap B=C\cap B=C$ is connected. 
 
 2) From 1) we have that $C\subset B$. If $B$ is  a point, then we must have $C=B$, which means that $\{x_n\}$ converges to the point $B$. 
 
 3) Since $A$ is at most countable, any connected component of $A$ is 1 point, and hence must be compact. If $C\not= \emptyset$, then $C$ must intersect at least one of the points in $A$, and hence by part 2) the whole sequence $\{x_n\}$ must converge to that point. Otherwise, $C=\emptyset$, and hence in this case $\{x_n\}$ diverges to infinity by definition.

\end{proof}

The following is used to prove properties of the Riemannian New Q-Newton's method algorithm. The notations $pr_{\pm}$ are from the Algorithm \ref{table:alg}.  

\begin{theorem} Let $X$ be a Riemannian manifold. Let $f:X\rightarrow \mathbb{R}$ be a $C^2$ function. Define $U=X\backslash \{x\in X:~grad(f)(x)=0\}$. For $x\in U$ and $\delta \in \mathbb{R}$ we define $A(x,\delta )=Hess(f)(x)+\delta ||grad (f)(x)||^{\alpha} Id$.

1) There exists a set $\Delta \subset \mathbb{R}$ of Lebesgue measure $0$ such that:  For all $\delta \in \mathbb{R}\backslash \Delta$, the set $\{x\in U:~A(x,\delta )$ is not invertible $\}$ has Lebesgue measure $0$. 

2) Fix $\delta \in \mathbb{R}\backslash \Delta$. We define $v(x,\delta )=A(x,\delta )^{-1}.grad (f)(x)$ and $w(x,\delta )=pr_{+,A(x,\delta )}.v(x,\delta )-pr_{-,A(x,\delta )}v(x,\delta )$, and $\lambda (x,\delta )=r(x)/||v(x,\delta )||$. Let $U_{\delta}=U\backslash \{x\in U:~\det (A(x,\delta ))=0\}$.  Then $w(x,\delta ),\lambda (x,\delta )$ are continuous on $U_{\delta}$. If we assume moreover that $f$ is $C^3$, then $w(x,\delta )$ is $C^1$ on $U_{\delta}$, and $\lambda (x,\delta )$ is  $C^1$ on $U_{\delta}\backslash \{x:~r(x)=||w(x,\delta )||\}$. 

3) Let the setting be as in part 2). Moreover, assume that $R$ is real analytic - like. Define $H(x,\delta )=R_x(-\lambda (x,\delta )w(x,\delta ))$, here $\lambda (x,s)$ is as in the definition for Riemannian New Q-Newton's method. There is a set $\Delta '\subset \mathbb{R}\backslash \Delta $ of Lebesgue measure $0$ so that for all $\delta \in \mathbb{R}\backslash (\Delta \cup \Delta ')$, the set $\{x\in U_{\delta }:~r(x)\not= ||w(x,\delta )||,~grad(H)(x,\delta )$ is not invertible$\}$ has Lebesgue measure $0$. 

4) Let the setting be as in part 3). There is a set $\Delta "\subset \mathbb{R}$ of Lebesgue's measure zero so that the following is satisfied. If $\delta \in \mathbb{R}\backslash \Delta "$ and $\mathcal{E}\subset X$ has Lebesgue measure $0$, then $H(.,\delta )^{-1}(\mathcal{E})\subset X$ also has Lebesgue's measure $0$.

\label{TheoremRandomnessLambda}\end{theorem}

\begin{proof}

1) Define $p(x,\delta )=\det (A(x,\delta ))$.  For $x\in U$, then $p(x,\delta )$ is a polynomial in $\delta$ of degree exactly $k=\dim (X)$. Note that $A(x,\delta )$ is not invertible if and only if $p(x,\delta )\not= 0$. We consider the set $\Gamma =\{(x,\delta )\in U\times \mathbb{R}:~p(x,\delta )=0\}$. Then $\Gamma $ is a measurable subset of $U\times \mathbb{R}$, and we need to show the existence of a set $\Delta \subset \mathbb{R}$ of Lebesgue measure $0$ so that for all $\delta \in \mathbb{R}\backslash \Delta$, the set $\Gamma _{\delta}=\{x\in U:~p(x,\delta )=0\}$ has Lebesgue measure $0$. 

First, we will show that $\Gamma$ has Lebesgue measure $0$. Let $1_{\Gamma}:~U\times \mathbb{R}\rightarrow \{0,1\}$ be the characteristic function of $\Gamma$.  For a set $A$, we denote by $|A|$ its Lebesgue measure. Then, by Fubini-Tonneli's theorem in integration theory for non-negative functions, one obtains: 
\begin{eqnarray*}
|\Gamma |=\int _{U\times \mathbb{R}}1_{\Gamma}(x,\delta ) d(x,\delta ) =\int _{U}(\int _{\mathbb{R} }1_{\Gamma}(x,\delta )d \delta )dx=0.
\end{eqnarray*}
This is because for each $x\in U$, the set $\{\delta \in \mathbb{R}:~p(x,\delta )=0\}\subset \mathbb{R}$ has cardinality at most $k$, and hence has Lebesgue measure $0$. 

Applying Fubini-Tonneli's theorem again, but now doing the integration on $U$ first, one obtains that for a.e. $\delta \in \mathbb{R}$, the set $\Gamma _{\delta}=\{x\in U:~p(x,\delta )=0\}\subset U$ has Lebesgue measure $0$.  

2) The claims follow from perturbation theory for linear operators \cite{kato}, where $pr_{\pm ,A(x,\delta )}$ can be represented via Cauchy's integrations on the complex plane containing the variable $\delta$. The readers can see for example \cite{truong-etal} on details how the arguments go. 

3)  For each $x\in U$, define $\Delta _x=\{\delta \in \mathbb{R}:~p(x,\delta )=0\}$. Then, $\Delta _x$ is a finite set (zeros of the polynomial $p(x,\delta )$).

By calculating, we see that when it is legit (c.f. part 2 above), then $grad(H(.,\delta ))$ is invertible at $x$ iff $q(x,\delta )=det (grad (H(.,\delta )))\not= 0$. 

Now, we fix $x_0\in U$ and $\lambda _0\in \Delta _{x_0}$. We note that there is an open interval $(a,b)\subset \mathbb{R}$ containing $\delta _0$, on which $grad (H(x,\delta ))$  is legit and $-\lambda (x,\delta )w(x,\delta )$ is a rational function of $\delta$, for $x$ close enough to $x_0$. Here is how to see this. We can write $\delta =\delta _0+\epsilon$, where $\epsilon$ is a small real number. Then, for $x$ close to $x_0$ and $\delta$ close to $\delta _0$, the eigenvalues of $A(x,\delta )$ are close to those of $A(x_0,\delta _0)$. Therefore, we can choose two close curves $\gamma ^+$ and $\gamma ^-$ in the complex plane containing the variable $\epsilon$, so that for $\delta =\delta _0+\epsilon$ close to $\delta _0$
\begin{eqnarray*}
w(x,\delta _0+\epsilon )&=& \frac{1}{2\pi i}\int _{\gamma ^+} ((A(x,\delta _0)+\epsilon ||grad(f)(x) ||^{\alpha})^{-1}-\zeta )^{-1}grad(f)(x) d\zeta \\
&&-\frac{1}{2\pi i}\int _{\gamma ^-} ((A(x,\delta _0)+\epsilon ||grad(f)(x) ||^{\alpha})^{-1}-\zeta )^{-1} grad(f)(x)d\zeta .
\end{eqnarray*}
 
Since we assumed that $||w(x_0,\delta _0)||\not= r(x_0)$, it follows that $||w(x,\delta _0+\epsilon )||\not= r(x)$ under the current assumptions. So we have that $\lambda (x,\delta _0+\epsilon )$ has the same form in the considered domain. Thus we can assume that $\lambda (x,\delta _0+\epsilon )=\gamma _{j+1}$ for a constant $j$. Moreover, the above formula for $w(x,\delta _0+\epsilon )$ clearly shows that it is a rational function on $\epsilon$, of degree exactly $m$. 

Hence, by the assumption that $R$ is real analytic-like, it follows that (see the discussion in Example \ref{Example4}), either $q(x_0,\delta _0+\epsilon )$ is zero identically, or its set of zero has Lebesgue's measure zero. We will show that the former case cannot happen. Indeed, while we defined the function 
\begin{eqnarray*}
\epsilon \mapsto \det (grad (R(-\lambda _jw(x,\delta _0+\epsilon ))))
\end{eqnarray*}
only for $\epsilon $ small, it is no problem to extend its domain of definition to the maximum domain  $W\subset \mathbb{R}$ for which $-\lambda _jw(x,\delta _0+\epsilon )\in B(x,r(x)) $. The same function $h_x\circ \kappa _x(s)$ (see Definition \ref{Def5}) will work for the whole domain $W$. 

Now, we note that for $\epsilon $ large enough (uniformly in $x$), then $-\lambda _{j+1}w(x,\delta _0+\epsilon ) \in B(x,r(x))$. It can be seen from the above integral representation for $w(x,s)$. Here is another, easier to see, way. We check this claim for example at the point $x_0$. Let  $\zeta _1,\ldots ,\zeta _l >0$ and $\zeta _{l+1},\ldots ,\zeta _{k}<0$ be all eigenvalues of $A(x_0,\delta _0)$. Let also $e_1,\ldots ,e_l$ and $e_{l+1},\ldots ,e_k$ be the corresponding orthonormal basis. If $grad(f)(x_0)=a_1e_1+\ldots +a_ke_k$, then 
\begin{eqnarray*}
w(x_0,\delta _0+\epsilon )=\sum _{\beta <l+1}\frac{a_{\beta}}{\zeta _{\beta } +\epsilon ||grad(f)(x_0)||^{\alpha }}e_{\beta}- \sum _{\beta >l}\frac{a_{\beta}}{\zeta _{\beta } +\epsilon ||grad(f)(x_0)||^{\alpha }}e_{\beta}. 
\end{eqnarray*}
Hence $\lim _{\epsilon\rightarrow\infty}-\lambda _{j+1}w(x_0,\delta _0+\epsilon )=0$, which confirms the claim.

Now, we are ready to show that $q(x,\delta _0+\epsilon )$ is not zero identically, for $\epsilon $ small enough. It suffices to consider the case $x=x_0$. Indeed, since $h_{x_0}\circ \kappa _{x_0}(\delta _0+\epsilon )$ restricts to $q(x_0,\delta _0+\epsilon )$, if $q(x_0,\delta _0+\epsilon )$ is zero identically, then the same would be true for $h_{x_0}\circ \kappa _{x_0}(\delta _0+\epsilon )$. However,  we will show that this is not the case. For $\epsilon $ large enough, from what written, we get that: $h_x\circ \kappa _x(\delta _0+\epsilon )=grad _x(R_x(-\lambda _{j+1}w(x,\delta _0+\epsilon )))$. On the other hand, by Chain's rule: 
\begin{eqnarray*}
grad _x(R_x(-\lambda _{j+1}w(x,\delta _0+\epsilon )))=(grad _xR_x|_{-\lambda _{j+1}w(x,\delta _0+\epsilon ))}).grad _x(-\lambda _{j+1}w(x,\delta _0+\epsilon ))).
\end{eqnarray*}
By explicit calculation, one then sees that $\lim _{\epsilon\rightarrow\infty}grad _x(R_x(-\lambda _{j+1}w(x,\delta _0+\epsilon )))=Id$. (One can check this readily for the Euclidean case, in which case $H(x,\lambda )=x-\lambda _{j+1}w(x,\delta  )$, using either representation of $w(x,\delta )$ above. The case of general Riemannian manifolds follows by choosing an appropriate coordinate chart.) Hence, in particular $\lim _{\epsilon \rightarrow \infty}h_{x_0} \circ \kappa _{x_0}(\delta _0+\epsilon )=1$, which shows that $h_{x_0}\circ \kappa _{x_0}(\delta _0+\epsilon )\not= 0$, as desired. 

From the representation of $w(x_0,\delta _0+\epsilon )$ above, we see that $$||w(x_0,\delta _0+\epsilon )||^2=\sum _{\beta }\frac{a_{\beta }^2}{(\zeta _{\beta }+\epsilon ||grad (f)(x_0)||^{\alpha})^2},$$ and hence clearly is not a constant. (For example, one can see that it is a rational function in $\epsilon _0$ with at least one pole, but is not identically $\infty$, because the limit when $\epsilon$ goes to $\infty$ is $0$.)  The same is true for all $x$ close enough to $x_0$. Therefore, for each $x$, the set of $\delta \in \mathbb{R}$ for which $||w(x,s)||=r(x)$ is at most a finite set. 

Combining all the above, here is a summary of what we proved so far: Except for a countable set of values $\delta \in \mathbb{R}$, the map $x\mapsto R_x(-\delta (x,\delta )v(x,\delta ))$ is $C^1$ and is a local diffeomorphism.  With this, now we finish the proof of this part 3). As in the proof of part 1), we define the set $\Gamma '$ to consist of all pairs $(x,\delta )$ for which $A(x,\delta )$ is invertible but the gradient of $H(x,\delta )$ is not invertible. The slices $\Gamma '_{x}=\{\delta :~(x,\delta )\in \Gamma '\}$ have been shown above to have Lebesgue's measure zero. Hence, by the use of Fubini-Tonelli's theorem, we get that the slices $\Gamma '_{\delta}=\{x:~(x,\delta )\in \Gamma '\}$ also have Lebesgue's measure zero for $\delta \in \mathbb{R}$ outside of a set of Lebesgue's measure $0$. 

 4) This is a consequence of part 3 and the following well-known fact: 
 
 {\bf Fact.} Let $H:U\rightarrow V$ be a map, where $U,V$ are open subsets in $\mathbb{R}^k$. Assume that there is a set $\mathcal{E}'\subset U$ of Lebesgue's measure $0$ so that for all $x\in U\backslash \mathcal{E}'$ then $H$ is $C^1$ near $x$ and $grad (H)$ is invertible there. Then, if $\mathcal{E}\subset V$ has Lebesgue's measure $0$, it follows that $H^{-1}(\mathcal{E})\subset U$ has Lebesgue's measure $0$.  

\end{proof}

The following two simple lemmas are useful for our purposes. 

\begin{lemma}
Let $U\subset \mathbb{R}^k$ be an open and convex set, and let $f:U\rightarrow \mathbb{R}$ be a $C^1$-function, whose gradient $\nabla f$ is Lipschitz continuous on $U$ with Lipschitz constant $L$. Let $x_0\in U$ and $\delta >0$ so that $x_0-\delta \nabla f(x_0)\in U$. Then 
\begin{eqnarray*}
f(x_0-\delta \nabla f(x_0))-f(x_0)\leq -\delta (1-\delta L) ||\nabla f(x_0)||^2. 
\end{eqnarray*}
\label{LemmaSimpleEstimate}\end{lemma}
\begin{proof}
By the Fundamental Theorem of Calculus, we have
\begin{eqnarray*}
f(x_0-\delta \nabla f(x_0))-f(x_0)=-\delta \int _0^1\nabla f (x_0-s\delta \nabla f(x_0)).\nabla f(x_0)ds. 
\end{eqnarray*}
Plugging into the RHS the following estimate
\begin{eqnarray*}
\nabla f (x_0-s\delta \nabla f(x_0)).\nabla f(x_0)&=&||\nabla f(x_0)||^2+(\nabla f (x_0-s\delta \nabla f(x_0)-\nabla f(x_0)).\nabla f(x_0))\\
&\geq & ||\nabla f(x_0)||^2-s\delta L||\nabla f(x_0)||^2,
\end{eqnarray*}
we obtain the result. 
\end{proof}

\begin{lemma}
Let $U\subset \mathbb{R}^k$ be an open and convex subset, and $H:U\rightarrow \mathbb{R}^k$ a continuous function. Assume that there is $\lambda >0$ so that $||H(x)-H(y)||\geq \lambda ||x-y||$ for all $x,y\in U$. Let $\mathcal{E}\subset \mathbb{R}^k$ be a measurable set of Lebesgue measure $0$. Then $H^{-1}(\mathcal{E})$ is also of Lebesgue measure $0$.
\label{LemmaMeasureEstimate}\end{lemma}
\begin{proof}
Since $\mathcal{E}$ has Lebesgue measure $0$, for each $\epsilon >0$ there is a sequence of balls $\{B(x_i,r_i)\}_{i=1,2,\ldots }$ so that $\mathcal{E}\subset \bigcup _{i=1}^{\infty}B(x_i,r_i)$ and $\sum _i(Vol(B(x_i,r_i)))\sim \sum _ir_i^k<\epsilon$. 

Since $||H(x)-H(y)||\geq \lambda ||x-y||$ for all $x,y\in U$, we get that $H^{-1}(\mathcal{E})$ $\subset$  $\bigcup _iB(H^{-1}(x_i),r_i/\lambda )$ and
\begin{eqnarray*}
\sum _iVol(B(H^{-1}(x_i),r_i/\delta ))\sim \sum _i (r_i/\delta )^k <\epsilon /\lambda ^k.
\end{eqnarray*} 
Therefore, the Lebesgue measure of $H^{-1}(\mathcal{E})$ is $<\epsilon /\lambda ^k$ for all $\epsilon >0$, and hence must be $0$. 
\end{proof}

\subsection{Global convergence guarantee to local minima for the new algorithms} Now we are ready to state and prove convergence and/or avoidance of saddle points for the mentioned algorithms.

We start with a global convergence guarantee for Riemannian Backtracking GD algorithm. 

\begin{theorem} Let $X$ be a Riemannian manifold and $r,R$ a Strong local retraction on $X$. Let $f:X\rightarrow \mathbb{R}$ be a $C^1$ function. Let $x_0\in X$, and $\{x_n\}$ the sequence constructed by the (Local) Riemannian Backtracking GD algorithm. Then,

1) If $x_{\infty}\in X$ is a cluster point of $\{x_n\}$, then $grad(f)(x_{\infty})=0$. 

2) Either $\lim _{n\rightarrow\infty}f(x_n)=-\infty$, or $\lim _{n\rightarrow\infty}d_X(x_n,x_{n+1})=0$. 

3) Let $A$ be the set of critical points of $f$, and let $B\subset A$ be a compact connected component. Let $C$ be the set of cluster points of $\{x_n\}$. If $C\cap B\not= \emptyset$, then $C\subset B$ and $C$ is connected. 

4) (Capture theorem.) Let the setting be as in part 3. If $B$ is a point, then $C=B$.

5) (Morse's function.) Assume that $f$ is a Morse function. Then the sequence $\{x_n\}$ either converges to a critical point of $f$, or diverges to infinity. Moreover, if $f$ has compact sublevels, then the sequence $\{x_n\}$ converges.

\label{Theorem11}\end{theorem} 
\begin{proof}

1) As in the Euclidean case (see e.g. \cite{bertsekas, truong-nguyen, truong-nguyen2}), by using Armijo's condition and Lemma \ref{LemmaLowerBound}, the proof follows. Here, one subtlety to note is that, by the properties of a Strong local retraction, if the sequence $\{x_{n_j}\}$ converges to $x_{\infty}$, then $||-\delta _{n_j}grad(f)(x_{n_j})||<3r(x_{\infty})/4$ when $j$ is large enough. Therefore, we can apply Armijo's condition for $f(x_{n_j+1})-f(x_{n_j})$ to obtain a contradiction if $\limsup _{j\rightarrow\infty}||grad(f)(x_{n_j})||$ were not zero.  

2) We use here the fact that $d_X(x_n,x_{n+1})$ is bounded by $||\delta _n grad(f)(x_n)||$, and proceed as in the Euclidean case. Note that here $\delta _n$ is also allowed to not be uniformly bounded,  see \cite{truong-nguyen, truong-nguyen2}. 

3) If $B$ is non-empty, then $\lim _{n\rightarrow\infty}f(x_n)>-\infty$. Hence, by 2) we have\\ $\lim _{n\rightarrow\infty}d_X(x_n,x_{n+1})=0$. Thus, Theorem \ref{TheoremConvergenceSequenceRiemannianManifolds} applies to give the desired conclusion. 

4) and 5) Follows from 3), as in the proof of Theorem \ref{TheoremConvergenceSequenceRiemannianManifolds}.

\end{proof}

In the literature, usually Capture theorem is stated under the assumption that $B$ is an isolated local minimum. Here we see that the isolated assumption alone is sufficient.

We next prove both global convergence guarantee and avoidance of saddle points for Riemannian Local Backtracking GD, and also for a continuous version for Riemannian Backtracking GD.  

\begin{theorem}
Let $f:X\rightarrow \mathbb{R}$ be a $C^{1}$ function, for which $\nabla f$ is locally Lipschitz continuous. Assume moreover that $f$ is $C^2$ near its generalised saddle points. 

A)

Then there is a smooth function $h:X\rightarrow (0,\delta _0]$ so that the map $H:X\rightarrow \mathbb{R}^k$ defined by $H(x)=R_x(-h(x)\nabla f(x))$ has the following property: 

(i) For all $x\in X$, we have $f(H(x)))-f(x)\leq -\alpha h(x)||\nabla f(x)||^2$. 

(ii) For every $x_0\in X$, the sequence $x_{n+1}=H(x_n)$ either satisfies $\lim _{n\rightarrow\infty}d_X(x_n,x_{n+1})=0$ or diverges to infinity. Each cluster point of $\{x_n\}$ is a critical point of $f$. a) If moreover, $f$ has at most countably many critical points, then $\{x_n\}$ either converges to a critical point of $f$ or diverges to infinity. b) More generally, if $B$ is a compact component of the set of critical points of $f$ and $C$ is the set of cluster points of $\{x_n\}$ and $B\cap C\not= \empty$, then $C\subset B$ and $C$ is a connected set. 

(iii) There is a set $\mathcal{E}_1\subset X$ of Lebesgue measure $0$ so that for all $x_0\in X\backslash \mathcal{E}_1$, if the sequence $x_{n+1}=H(x_n)$ {\bf converges}, then the limit cannot be a {\bf generalised} saddle point.  

(iv) There is a set $\mathcal{E}_2\subset X$ of Lebesgue measure $0$ so that for all $x_0\in X\backslash \mathcal{E}_2$, any cluster point of the sequence $x_{n+1}=H(x_n)$ is not a saddle point, and more generally cannot be an isolated generalised saddle point. 

B) The statements i), ii), iii) and iv) above also hold for Riemannian Local Backtracking GD.

\label{Theorem3}\end{theorem}
\begin{proof}

A) For the continuous version of Riemannian Backtracking GD: 

{\bf The Euclidean setting $X=\mathbf{R}^k$:} We will first present the proof for the Euclidean case, which is simpler, and then indicate at the end what needs to be modified in the general Riemanian manifold setting.

Since $\nabla f$ is locally Lipschitz continuous, for each $x\in \mathbb{R}^k$, there are positive numbers $r(x),L(x)>0$ so that  $\nabla f(x)$ has Lipschitz constant $L(x)$ in the ball $B(x,r(x))$. That is, for all $y,z\in B(x,r(x))$ we have $||\nabla f(y)-\nabla f(z)||\leq L(x)||y-z||$. Also, we can choose $L(x)$ large enough so that for all $0<\delta \leq 1/(L(x))$ then $f(z-\delta \nabla f(z))-f(z)\leq -\alpha \delta ||\nabla f(z)||^2$ for all $z\in B(x,r(x))$.  

 There is a partition of unity $\{\varphi _j\}_{j=1,2,\ldots }$ of $\mathbb{R}^k$ with compact supports so that for every $j\in \mathbb{N}$, there is a point $z_j\in \mathbb{R}^k$ for which the support $supp(\varphi _j)$ of $\varphi _j$ is contained in $B(z_j,r(z_j))$. Moreover, $\{supp (\varphi _j)\}_j$  is locally finite, that is  every $x\in \mathbb{R}^k$ has an open neighbourhood $U$ which intersects only a finite number of those $supp (\varphi _j)$'s. (This is related to Lindel\"off theorem. We recall here the main idea: Any open cover of $\mathbb{R}^k$ has a subcover which is locally finite. We can even arrange that each point $x\in \mathbb{R}^k$ is contained in at most $k+1$ open sets in the subcover. Then we construct a partition of unity with compact supports contained in open sets in that subcover.). For each $j=1,2,\ldots $, we let 
\begin{eqnarray*}
M_j=\max _{y_1,y_2\in B(z_j,r(z_j))}||\nabla \varphi _j(y_1)|| \times ||\nabla f(y_2)||. 
\end{eqnarray*}  
  
 We now define the function $h:\mathbb{R}^k\rightarrow \mathbb{R}$ by the following formula: 
\begin{eqnarray*}
h(x):=\sum _{j=1}^{\infty} \frac{1}{10^{j}(M_j+1)} \varphi _j(x)\min \{\frac{1}{2L(z_j)},\delta _0,1\}.
\end{eqnarray*}

Since $\varphi _j$'s are non-negative, smooth and their supports are locally finite, it follows that the function $h(x)$ is well-defined, smooth and non-negative. Since $\sum _j\varphi _j(x)=1$, it follows that $0<h(x)\leq \delta _0$ for all $x\in \mathbb{R}^k$. 

Fix a point $x\in \mathbb{R}^k$. Then there is at least one $j$ so that $\varphi _j(x)>0$, and hence $h(x)\geq \frac{1}{10^j(M_j+1)} \varphi _j(x)\min \{\frac{1}{2L(z_j)},\delta _0,1\} >0$.  Then there is a finite index set $J$ and an open neighbourhood $V$ of $x$ so that if $supp (\varphi _j)\cap V\not= \emptyset$ then $j\in J$. Since $supp(\varphi _j)$'s are all compact, we can shrink $V$ so that $x\in \bigcap _{i\in J}B(z_j,r(z_j))$ and $h(x)\leq \max _{j\in J}\min \{\frac{1}{2L(z_j)},\delta _0\}$. Since $h$ is a smooth function, we can find an open neighbourhood $U$ of $x$ so that $U\subset  V\cap \bigcap _{i\in J}B(z_j,r(z_j))$ and for all $y\in U$ we have $h(y)\leq  \max _{j\in J}\min \{\frac{2}{3L(z_j)},\delta _0\}$. It then follows by the choice of $L(z_j)$'s that for all $y\in U$ we have
\begin{eqnarray*}
f(y-h(y)\nabla f(y))-f(y)\leq -\alpha h(y) ||\nabla f(y)||^2. 
\end{eqnarray*}   
We will now show also that $H(y)=y-h(y)\nabla f(y)$ is injective on this same set $U$. In fact, assume that there are distinct $y_1,y_2\in U$ so that  $y_1-h(y_1)\nabla f(y_1)=y_2-h(y_2)\nabla f(y_2)$. We rewrite this as: 
\begin{eqnarray*}
(y_1-y_2)-h(y_1)(\nabla f(y_1)-\nabla f(y_2))-(h(y_1)-h(y_2))\nabla f(y_2)=0. 
\end{eqnarray*}
Since $U\subset V$ and by the choice of $V$, if $\varphi _j(y)\not= 0$ for some $y\in V$ then $j\in J$.  Hence, we obtain a contradiction because
\begin{eqnarray*}
||h(y_1)(\nabla f(y_1)-\nabla f(y_2))||\leq \frac{2}{3}||y_1-y_2||,
\end{eqnarray*}
and using that $|\varphi _j(y_1)-\varphi _j(y_2)|\leq ||y_1-y_2|| \times \max _{B(z_j,r(z_j))}||\nabla \varphi _j||$
\begin{eqnarray*}
&&||(h(y_1)-h(y_2))\nabla f(y_2)||\\
&=&|\sum _{j\in J}  \frac{1}{10^{j}(M_j+1)} (\varphi _j(y_1)-\varphi _j(y_2))\min \{\frac{1}{2L(z_j)},\delta _0,1\}|\times ||\nabla f(y_2)||\\
&\leq&\sum _{j\in J}  \frac{1}{10^{j}(M_j+1)} M_j\min \{\frac{1}{2L(z_j)},\delta _0,1\}||y_1-y_2||\\
&\leq& \frac{1}{9}||y_1-y_2||. 
\end{eqnarray*}

Near generalised saddle points of  $f$, the map $x\mapsto H(x)$ is moreover $C^1$ by the assumption that $f$ is $C^2$ near those points. Hence, $x\mapsto H(x)$ is a local diffeomorphism near generalised saddle points of $f$. 

{\bf Proof of (i)}: Already given above. 

{\bf Proof of (ii)}: Since $h(x)$ is smooth and positive, it follows that for every compact subset $K\subset \mathbb{R}^k$, we have $\inf _{x\in K}h(x)>0$. Hence, we can use the same arguments as in \cite{bertsekas} to show that any cluster point of the sequence $\{x_{n+1}=H(x_n)\}$ is a critical point of $f$. 

Since $h(x)\leq \delta _0$ for all $x\in \mathbb{R}^k$, we can use the same proof of part 1 of Theorem 2.1 in \cite{truong-nguyen} to conclude that either $\lim _{n\rightarrow\infty}||x_{n+1}-x_n||=0$ or $\lim _{n\rightarrow\infty}||x_n||=\infty$. 

Then we can use the same proof of part 2 of Theorem 2.1 in \cite{truong-nguyen}, by employing the real projective space $\mathbb{P}^k$ and result in \cite{asic-adamovic} for cluster points of sequences $\{x_n\}$ in compact metric space $(X,d)$ satisfying $\lim _{n\rightarrow\infty}d(x_{n+1},x_n)=0$, to show that if $f$ has at most countably many critical points, then either $\{x_n\}$ converges to a critical point of $f$ or $\lim _{n\rightarrow\infty}||x_n||=\infty$. 

{\bf Proof of (iii)}: By using Stable-Center Manifold theorem for local diffeomorphisms and using that the map $x\mapsto H(x)$ is a local diffeomorphism near generalised saddle points, we can argue as in \cite{lee-simchowitz-jordan-recht, panageas-piliouras} that there is an open neighbourhood $U$ of generalised saddle points of $f$, and a subset $\mathcal{F}_1\subset U$ of Lebesgue measure $0$, so that if all $\{x_n\}\subset U\backslash \mathcal{F}_1$ and $\{x_n\}$ {\bf converges}, then the limit point cannot be a saddle point.   

Since we showed in the construction that the map $x\mapsto H(x)$ is locally injective, it follows that $\mathcal{E}_1=\bigcup _{n\in \mathbb{N}}H^{-n}(\mathcal{F}_1)$ also has Lebesgue measure $0$. Then it follows that if $x_0\in \mathbb{R}^k\backslash \mathcal{E}_1$, then $\{x_n\}$ cannot converge to a generalised saddle point. 

{\bf Proof of (iv)}: We use the ideas in \cite{truong-nguyen}. Note that a saddle point of $f$ is a non-degenerate critical point, and hence is an isolated generalised saddle point. Assume that the set of cluster points $A$ of $\{x_n\}$ contains an isolated generalised saddle point $y_0$. Then the property $\lim _{n\rightarrow\infty}||x_n||=\infty$ does not hold, and hence by part (ii) we must have $\lim _{n\rightarrow\infty}||x_{n+1}-x_n||=0$. Then, by the result in \cite{asic-adamovic} for the real projective space $\mathbb{P}^k$, it follows that the closure $\overline{A}\subset \mathbb{P}^k$ is connected. By part (ii) again, the set A is contained in the set of all critical points of $f$, and hence $y_0$ is also an isolated point of A, and hence of $\overline{A}$. Thus $A=\{y_0\}$, and hence $\lim _{n\rightarrow\infty}x_n=y_0$. Then we can use part (iii) to conclude.   

{\bf The general Riemannian manifold setting:} The proof is similar to the above Euclidean setting. Here are some necessary modifications. 

The main point is to use the stronger form of Lindel\"off lemma that, since $X$ is a Riemannian manifold of finite dimension $k$, each open covering of $X$ has a subcovering which is locally finite. The latter means that for every $x$, there is a small open subset $U$ around $x$, so that every point $y\in U$ belongs to at most $k+1$ sets in the subcovering.  We then use this fact and partition of unity to cook up a smooth function for learning rates $h:X\rightarrow (0,\delta _0]$ which has the same properties as that of the learning rates in the proof of Theorem \ref{Theorem3}. When having this, we can proceed as before. 

B) For the Riemannian Local Backtracking GD:

Note that by construction we have $\widehat{\delta} (x_n)\leq \delta _0$ for all $n$. Since $L(x), r(x)$ and $\nabla f(x)$ are continuous in $x$, it follows that for all compact subset $K\subset \mathbb{R}^k$, we have $\inf _{x\in K}\delta (x)>0$. Note that by construction $x-\widehat{\delta}(x)\nabla f(x)\in B(x,r(x)/2)$ for all $x\in \mathbb{R}^k$ and hence Armijo's condition is satisfied. Therefore, i) and ii) follow as in the case for Backtracking GD on Euclidean spaces (\cite{truong-nguyen, truong-nguyen2}). 

Now we prove iv). By using Lindeloff's trick as in \cite{panageas-piliouras}, we need to work with a countable set $\mathcal{S}$ of saddle points of $f$. Therefore, for a {\bf random choice} of $\alpha, \beta, \delta _0$, we have the following: for every $x\in \mathcal{S}$ then $\alpha /L(x)$ does not belong to $\{\beta ^n\delta _0: ~n=0,1,2,\ldots \}$, and also $r(x)/||\nabla f(x)||$ does not belong to $\{\beta ^n\delta _0: ~n=0,1,2,\ldots \}$. 

Hence, for each $x\in \mathcal{S}$, either $\alpha /L(x)>\delta _0$, or there is a number $n(x_0)$ so that
\begin{eqnarray*}
\beta ^{n(x_0)+1}\delta _0<\frac{\alpha }{L(x)}< \beta ^{n(x_0)}\delta _0. 
\end{eqnarray*}
In both cases, since $z\mapsto L(z)$ is continuous, there is an open neighbourhood $U(x)$ of $x$ so that for all $z\in U(x)$, then $\alpha /L(z)$ has the same behaviour.  Shrinking $U(x)$ if necessary,  we can assure that $||\nabla f(z)||$ is small in $U(x)$, and hence $r(z)/||\nabla f(z)||>\delta _0$.Therefore, by Definition \ref{DefinitionBacktrackingGDNew}, we have that $\widehat{\delta} (z)=\widehat{\delta} (x)=\beta ^{n(x_0)+1}\delta _0$ for all $z\in U(x)$. In particular, the map $z\mapsto H(z)=z-\widehat{\delta} (z)\nabla f(z)$ is a local diffeomorphism in $U(x)$. 

By \cite{truong-nguyen}, if the cluster set of $\{x_n\}$ contains an isolated generalised saddle point, then $\{x_n\}$ converges to that generalised saddle point. Then, we can apply Stable-Central theorem to obtain that there is an open neighbourhood $U$ of $\mathcal{S}$ and a set $\mathcal{F}\subset U$ of Lebesgue measure $0$ such that if $x_0\in \mathbb{R}^k$ is so that the cluster points of $\{x_n\}$ contains an isolated generalised saddle point, then there is $n(x_0)$ for which $H^{n(x_0)}(x_0)\in \mathcal{F}$. 

We note that since $z\mapsto L(z)$, $z\mapsto r(z)$ and $z\mapsto ||\nabla f(z)||$ are continuous, for every $x\in \mathbb{R}^k$, there is a neighbourhood $U(x)$ so that for all $z\in U(x)$ then $\widehat{\delta} (z)=\widehat{\delta} (x)$ or $\widehat{\delta} (z)=\widehat{\delta} (x)/\beta$. Then we see that $z\mapsto H(z)$ is injective on both sets $\{z\in U(x): ~\widehat{\delta} (z)=\widehat{\delta} (x)\}$ and $\{z\in U(x): ~\widehat{\delta} (z)=\widehat{\delta} (x)/\beta \}$. Then argue as in the proof of (iii) and (iv) in A), we get that the set 
\begin{eqnarray*}
\mathcal{E}_2=\bigcup _{n\in \mathbb{N}}H^{-n}(\mathcal{F})
\end{eqnarray*}
has Lebesgue measure $0$, and for all $x_0\in \mathbb{R}^k\backslash \mathcal{E}_2$, the sequence $\{x_n\}$ constructed from Definition \ref{DefinitionBacktrackingGDNew} cannot have any cluster point which is an isolated generalised saddle point. 

The proof of iii) is similar.

\end{proof}

Finally, we state and prove the result for Riemannian Backtracking New Q-Newton's method. 

\begin{theorem} Let $f:X\rightarrow \mathbb{R}$ be a $C^3$ function. Let $\{x_n\}$ be a sequence constructed by the Riemannian Backtracking New Q-Newton's method. Then

1) Any cluster point $x_{\infty}$ of the sequence $\{x_n\}$ satisfies $grad (f)(x_{\infty})=0$, that is $x_{\infty}$ is a critical point of $f$.

2) Assume that $r,R$ satisfy the Real analytic-like condition.  There is a set $\mathcal{A}\subset X$ of Lebesgue measure $0$, so that if $x_0\notin \mathcal{A}$, then $x_{\infty}$ cannot be  a saddle point of $f$. 

3) Assume that $r,R$ satisfy the Real analytic-like condition. If $x_0\notin \mathcal{A}$ (as defined in part 2) and $Hess(f)(x_{\infty})$ is invertible, then $x_{\infty}$ is a local minimum and the rate of convergence is quadratic. 

4) More generally, if $Hess(f)(x_{\infty})$ is invertible (but no assumption on the randomness of $x_0$), then the rate of convergence is at least linear. 

5) If $x_{\infty}'$ is a non-degenerate local minimum of $f$, then for initial points $x_0'$ close enough to $x_{\infty}'$, the sequence $\{x_n'\}$  constructed by Riemannian New Q-Newton's method will converge to $x_{\infty}'$. 
\label{Theorem4}\end{theorem} 
\begin{proof}
The proof is similar to that of the Euclidean counterpart in \cite{truong1}. We note that in parts 2) and 3), if the Strong local retraction $r,R$ is not required to satisfy the Real analytic-like condition, then one still can prove that there are local Stable-Center manifolds for the associated dynamical systems near saddle points. Hence, on all Riemannian manifolds one has local guarantee for avoidance of saddle points near saddle points. The Real analytic-like condition is needed to assure global avoidance of saddle points, via the use of Theorem \ref{TheoremRandomnessLambda} Parts 4) and 5), which are local in nature, can be proven exactly as in \cite{truong1, truong-etal}. 

There is a subtle point in the proof of part 1), compared to its Euclidean counterpart, which lies in the fact that $\lambda _n $ in the update rule is in general not the constant $1$. This is because in general we do not have global retractions, for example if we work with open subsets of complete Riemannian manifolds. The choice of $\lambda _n$ in our update rule, which is about $1/||w_n||$ when $||v_n||$ is large, is important to assure that $grad (f)(x_{\infty})=0$. If $\lambda _n$ has another asymptotic growth, such as $1/||w_n||^2$, then there is no such guarantee.  Here we give a detailed proof of 1) to illustrate the point. 

Proof of part 1):  For simplicity, we assume that the whole sequence $\{x_n\}$ converges to $x_{\infty}$, so can adapt the proof in \cite{truong-etal} (and in the general case where we modify by following the proof in \cite{truong1}). Since $x_{n+1}=R_{x_n}(-\lambda _nw_n)$ converges to $x_{\infty}$, and $||\lambda _nw_n||\leq r(x_n)/2$ for all $n$, together with the fact that $R_{x_{\infty}}(v)$ is a diffeomorphism for $||v||<r(x_{\infty})$, it follows that $\lim _{n\rightarrow\infty}||-\lambda _nw_n||=0$. This implies that first of all, $||w_n|| $ is bounded, since if it were true that $\lim _{n\rightarrow\infty}||w_n||=\infty$, then for large $n$ we would have $\lambda _n\sim r(x_n)/(2||w_n||)$, and hence we would have a contradiction that $\lim _{n\rightarrow\infty}||\lambda _nw_n||>0$. Therefore, $\lambda _n$ is uniformly bounded from below by a positive number, thus from $\lim _{n\rightarrow\infty}||-\lambda _nw_n||=0$ we obtain also that $\lim _{n\rightarrow\infty}||w_n||=0$. Since $||v_n||=||w_n||$ for all $n$, we have that $\lim _{n\rightarrow\infty}||v_n||=0$. Then (note that $||A_n||$ is uniformly bounded) $$||grad(f)(x_{\infty})||=\lim _{n\rightarrow\infty}||grad(f)(x_n)||=\lim _{n\rightarrow\infty}||A_n.v_n||=0,$$ 
as wanted. 
\end{proof}

\subsection{Application: Solving optimisation problems on the unit sphere} In many cases, one has the need to solve optimisation/root-finding on the unit sphere. For example, this is the case where one wants to find eigenvectors of a square matrix. In this subsection, we propose a method for this, by applying the Riemannian versions of Backtracking line search.   

For optimisation problems of the form $\min _{||x||\leq 1}f(x)$, where $x\in \mathbb{R}^{k+1}$, methods in the Euclidean space may not be suitable. Here, we show that combining that with Riemnannian optimisation on the sphere $\{x:~||x||=1\}$ can yield improved performance. 

{\bf General strategy: }

Do Riemannian optimisation on the manifold $S^{k}=\{x\in \mathbb{R}^{k+1}:~||x||=1\}$, with $r(x)=\pi$. (Note that, in experiments, we see that even with putting $r(x)=\infty$, still the performance is very good.) Here, there are two ways to choose $R_x(v)$:

Way 1: $R_x(v)=(x+v)/\sqrt{1+||v||^2}$. 

Way 2: (geodesic) $R_x(v)=\cos (||v||)x+\sin (||v||)v/||v||$. 

We illustrate this with the problem of finding eigenvectors of a symmetric matrix with real coefficients.  It is equivalent to considering a quadratic function $f_A(x)=<Ax,x>/2$, where $A$ is a symmetric matrix. In this case, if $\lambda _1(A)$ is the smallest eigenvalue of $A$, then $\min _{x\in S^{k}}f(x)=\lambda _1(A)/2$. Hence, this problem is interesting also for numerical linear algebra. 

For $S^{k}$, we will use the induced metric from $\mathbb{R}^{k+1}$. This implies, in particular that if $v\in T_xS^{k}$, then $||v||_{T_xS^{k}}=||v||_{R^{k+1}}$.  The computations for Riemannian gradient and Hessian are also quite nice, the next 2  formulas are taken from \cite[Propositions 3.49, Section 5.5]{boumal}: If $x\in S^{k}$ and $v\in T_xS^{k}$, then (The RHS of the formulas are interpreted in the usual Euclidean setting)
\begin{eqnarray*}
grad(f_A)(x)&=&Ax-<Ax,x>x,\\
Hess(f_A)(x)[v]&=&Av-<Av,x>x-<Ax,x>v.
\end{eqnarray*} 
 
 As we mentioned above, the Riemannian Hessian is symmetric on $T_xS^{k}$. On the other hand, its obvious extension (using the same formula on the RHS) to $T_x\mathbb{R}^{k+1}$ may not be symmetric. Since it would be more convenient to do calculations with a {\bf symmetric} extension of $Hess(f_A)(x)$ to the whole $T_x\mathbb{R}^{k+1}$ (for example, when we want to decompose into positive and negative eigenvalues as in New Q-Newton's method), we will define explicitly such an extension $B:T_x\mathbb{R}^{k+1}\rightarrow T_x\mathbb{R}^{k+1}$. The most convenient way is to use, for $v\in T_x\mathbb{R}^{k+1}$, its orthogonal projection $v-<v,x>x$ to $T_xS^{k}$. Hence, we define $B$ by the formula:
 \begin{eqnarray*}
 B[v]:=Hess(f_A)(x)[v-<v,x>x]. 
 \end{eqnarray*}
 
With the above formulas, we can apply Riemannian Backtracking GD and Riemannian New Q-Newton's method, the Riemannian Newton's method and its random damping version, as well as Riemann Standard GD. For the learning rate for Riemannian Standard GD, we fix it to be $0.001$. We see from the experiments that the convergence here is very fast. 

\begin{table}[htp]
\fontsize{9}{9}\selectfont
  \centering
  \begin{tabular}{|l|c|c|c|c|c|c|}
  \hline
~~~~~~~Method/Function	& Example 1&Example 2   & Example 3  \\
\hline

Riem Newton & 10/(0.70,0.70) & 10/(-0.13,0.73,0.66)/S& 10/(-0.13,0.73,0.66)/M \\
Riem New Q & 10/(-0.70,0.70)/M &10/(-0.33,-0.66,0.66)/M & 10/(-0.13,0.73,0.66)/M \\
Riem Rand Newton & 10/(0.70,0.70)& 10/(-0.13,0.73,0.66)/S & 10/(-0.13,0.73,0.66)/M \\
Riem Backtrack &3/(-0.70,0.70)/M  &50/(-0.33,-0.66,0.66)/M &10/(-0.13,0.73,0.66)/M \\
Riem SGD &10/(0.42,0.90)  &10/(0.32,-0.24,0.91) &10/(0.017,0.701,0.71) \\
\hline
  \end{tabular}
   \caption{Results of experiments on finding the minimum eigenvalue and a corresponding eigenvector of a square matrix, as a Riemannian optimisation problem on the {\bf unit sphere}, with its induced metric from the Euclidean metric. The maximum number of iterates is 10, but the algorithm can stop before that because either the size of the gradient is smaller than a threshold ($1e-10$), or there is an error. The format is n/x/Remarks, where n is the number of iterates needed to achieve the point $x$. Legends: "E" for errors, "D" for divergence, "C" for convergence, "Newton" for Newton's method, "New Q" for New Q-Newton's method,  "Backtrack" for Backtracking GD, "SGD" for Standard GD, "Rand" for random, "Riem" for Riemannian, "S" for (near ) a saddle point or local maximum, "M" for (near) a local minimum.}   
  \label{tab:Tasks3}
\end{table}

{\bf Example 1.} Consider $\min _{||x||=1}f_A(x)$, where $A$ is the matrix 

 \[ \left( \begin{array}{cc}
2&4\\
4&2\\
\end{array}\right) \]

 The initial point will be $x_0/||x_0||=(0.4472136  ,0.89442719)\in S^1$.

{\bf Example 2.} Consider $\min _{||x||=1}f_A(x)$, where $A$ is the matrix: 

 \[ \left( \begin{array}{ccc}
-23&-61&40\\
-61&-39.5&155\\
40&155&-50\\
\end{array}\right) \]

 The initial point will be $x_0/||x_0||=(0.29369586,$ $0.54091459,$ $0.78813333)\in S^2$.

{\bf Example 3:} Consider $\min _{||x||=1}-f_A(x)$, where $A$ is the matrix in Example 2. We start from the same initial point $x_0/||x_0||=(0.29369586,$ $0.54091459,$ $0.78813333)\in S^2$. Here, we want to find the largest eigenvalue of $A$, hence we consider $-f_A(x)$ and not $f_A(x)$.

\section{Optimisation/Root-finding on Banach spaces}

Here we aim to extend the above results to infinite dimensional Banach spaces. For a comprehensive reference on Banach spaces, we refer the readers to \cite{rudin}. We recall that a Banach space $X$ is a vector space, together with a norm $||.||$, and is complete with respect to the norm. The topology induced by $||.||$ is called the strong topology. The dual space $X^*$, consisting of bounded linear maps $h:X\rightarrow \mathbb{R}$, is also a Banach space, with norm $||h||=\sup _{x\in X:~||x||\leq 1}|h(x)|$.   

There are many pathologies when one works with infinite dimensional Banach spaces. For example, a bounded sequence $\{x_n\}$ may not have any subsequence converging in the (strong) topology, unlike the case of finite dimensions where we have Bolzano-Weierstrass property. For to better deal with this, we need to work with another topology, called weak topology. To distinguish with the usual (strong) topology, we will use a prefix "w-" when speaking about weak topology. Here, a sequence $\{x_n\}$ w-converges to $x$ if for every $h\in X^*$ we have $\lim _{n\rightarrow\infty}h(x_n)=h(x)$. Note that strong and weak topologies coincide only if $X$ is finite dimensional. If we want to preserve the good Bolzano-Weierstrass property, then we need to require that $X$ is reflexive, in that $X$ is isomorphic to its double dual $X^{**}$. In this case, it follows by Kakutani's theorem that any bounded sequence $\{x_n\}$ has a subsequence which is w-convergent. 

Finding critical points of functions in infinite Banach spaces has many applications, among them is finding solutions to PDE. For example, a fundamental technique in PDE is to use a weak formulation. For example, $u$ is a $C^2$ solution to Poisson's equation $-\nabla ^2u=g$ on a domain $\Omega \subset \mathbb{R}^d$ with $u|_{\partial \Omega }=0$, if and only if for all smooth functions $v$ with the same boundary condition we have
\begin{eqnarray*}
-\int _{\Omega}(\nabla ^2u)vdx = \int _{\Omega}gv.
\end{eqnarray*}
 Now, by using integration by parts, one reduce the above to 
 \begin{eqnarray*}
 \int _{\Omega}\nabla u\nabla vdx-\int _{\Omega}gvdx=0.
 \end{eqnarray*}
 One can then construct a Sobolev's space of measurable functions $H^{1}_0(\Omega)$ of functions with weak zero boundary conditions, and define a norm
 \begin{eqnarray*}
 ||u||_*^2=\int _{\Omega}||\nabla u(x)||^2dx.
 \end{eqnarray*}
It turns out that $H^{1}_0(\Omega)$ is a Banach space, and is indeed a Hilbert space, and hence is reflexive. To find a solution to the original Poisson's equation, one proceeds in  two steps. In Step 1, one finds a weak solution $u$ in $H^{1}_0(\Omega)$, which is in general only a measurable function. In Step 2, one uses regularity theory to show that solutions from Step 1 are indeed in $C^2$ and hence are usual solutions. Now, we can relate Step 1 to finding critical points in the following way. We define a function $f:H^{1}_0(\Omega )\rightarrow \mathbb{R}$ by the formula
\begin{eqnarray*}
f(u)=\frac{1}{2}\int _{\Omega}||\nabla u(x)||^2dx-\int _{\Omega}gu.
\end{eqnarray*}
 We can check that $f$ is a $C^1$ function, and that $u$ is a weak solution of Poisson's equation iff $u$ is a critical point of $f$. 
 
Because of their importance in many aspects of science and life, there have been a lot of research on finding solutions of PDE. We concentrate here the approach of using numerical methods to solve PDE. There is a comprehensive list of such methods, such as Finite difference method, Method of lines, Finite element method, Gradient discretization method and so on. It seems, however, that numerical methods for finding critical points, which are powerful tools in finite dimensional spaces, are so far not much used. One reason may well be that numerical methods for finding critical points themselves are not much developed in infinite dimension. Indeed, it seems that most of the (convergence) results in this direction so far is based on Banach's fixed point theorem for contracting mappings, and applies only to maps which have some types of monotonicity. For systematic development in GD, one can see for example the very recent work \cite{geiersbach-scarinci}, where the problem is of stochastic nature, the space is a Hilbert space, and $f$ is either in $C^{1,1}_L$ (in which case only the results $\lim _{n\rightarrow\infty}||x_{n+1}-x_n||=0$ is proven); or when with additional assumption on learning rates $\sum _n\delta _n=\infty$ and $\sum _{n}\delta _n^2<\infty$, where convergence to $0$ of $\{\nabla f(x_n)\}$ can be proven. Some other references are \cite{gallego-etal} (where the function is assumed to be strongly convex) and \cite{blank-rupprecht} (where the VMPT method is used), where again only $\lim _{n\rightarrow\infty}||\nabla f(z_n)||=0$ is considered, and no discussion of (weak or strong) convergence of $\{x_n\}$ itself or avoidance of saddle points are given. 

In this subsection, we will develop more the Backtracking GD method in infinite dimensional Banach spaces. To motivate the assumptions which will be imposed later on the functions $f$, we present here some further pathologies one faces in infinite dimensions. The first pathology is that for a continuous function $f:X\rightarrow \mathbb{R}$, it may happen that $f$ is unbounded on bounded subsets of $X$. The second pathology is that it may happen that $\{x_n\}$ w-converges to $x$, but $\nabla f(x_n)$ does not w-converge to $\nabla f(x)$. One needs to avoid this pathology, at least in the case when additionally $\lim _{n\rightarrow\infty}||\nabla f(x_n)||=0$, since in general the sequence $\{x_n\}$ one constructs from a GD method is not guaranteed to (strongly) converge, and hence at most can be hoped to w-converge to some $x$. If it turns out  that $\nabla f(x)\not= 0$, one cannot use GD method to find critical points. The second pathology is prevented by the following condition: 

{\bf Condition C.} Let $X$ be a Banach space. A $C^1$ function $f:X\rightarrow \mathbb{R}$ satisfies Condition C if whenever $\{x_n\}$ weakly converges to $x$ and $\lim _{n\rightarrow\infty}||\nabla f(x_n)||=0$, then $\nabla f(x)=0$.   

From a purely optimisation interest, Condition C is satisfied by quadratic functions, such as those considered on page 63 in \cite{nesterov} where Nesterov constructs some examples to illustrate the pathology with using GD in infinite dimensions. Indeed, if $f\in C^1$ is quadratic, then there is $A:X\rightarrow X^*$ a bounded linear operator and $b\in X^*$, so that $\nabla f(x)=Ax+b$. It is known that $A$ is then continuous in the weak topology. Therefore, if $\{x_n\}$ w-converges to $x$, then $\{Ax_n+b\}$ w-converges to $Ax+b$. In particular, $0=\liminf _{n\rightarrow\infty}||\nabla f(x_n)||\geq ||\nabla f(x)||$. Thus we have $\nabla f(x)=0$, and hence Condition C is satisfied.  

More interestingly, Condition C is valid for functions $f\in C^1$ which are convex. Indeed, assume that $f$ is in $C^1$ and is convex, and that $\{x_n\}$ w-converges to $x$ so that $\lim _{n\rightarrow\infty}||\nabla f(x_n)||=0$. Then since $f$ is continuous, it follows that $f$ is weakly lower semicontinuous (see Proposition 7b in \cite{browder}), that is $\liminf _{n\rightarrow\infty}f(x_n)\geq f(x)$. Now, since $f$ is convex and $\lim _{n\rightarrow\infty}||\nabla f(x_n)||=0$, for every $y\in X$ we have
\begin{eqnarray*}
f(y)\geq \limsup _{n\rightarrow\infty} [f(x_n)+<\nabla f(x_n),x_n-y>]=\limsup _{n\rightarrow\infty} f(x_n)\geq f(x).
\end{eqnarray*}
Therefore, $x$ is a  (global) minimum of $f$ and hence $\nabla f(x)=0$. 

It can be checked that the function associated to Poisson's equation satisfies Condition C (for example, because it is quadratic). More generally, the following class of non-linear PDE, see Definition 2 in \cite{browder}, satisfies Condition C:

{\bf Class $(S)_+$.} A $C^1$ function $f:X\rightarrow \mathbb{R}$ is of class $(S)_+$ if whenever $\{x_n\}$ w-converges to $x$ and $\limsup _{n\rightarrow\infty}<\nabla f(x_n),x_n-x>\leq 0$, then $\{x_n\}$ (strongly) converges to $x$. 

Indeed, assume that $f$ is in Class $(S)_+$, and $\{x_n\}$ w-converges to $x$ so that $\lim _{n\rightarrow\infty}$ $||\nabla f(x_n)||$ $=0$. Then, since $<\nabla f(x_n),x_n-x>\leq ||\nabla f(x_n)||\times ||x_n-x||$ and $\{||x_n-x||\}$ is bounded,  we have $\limsup _{n\rightarrow\infty}<\nabla f(x_n),x_n-x>\leq 0$. Therefore, by definition of Class $(S)_+$, we have $\{x_n\}$ (strongly) converges to $x$. Then, $\nabla f(x)=\lim _{n\rightarrow\infty}\nabla f(x_n)=0$. Hence Condition C is satisfied. ({\bf Remark.} Here we obtain the same conclusion under the weaker condition that $\nabla f$ is only demi-continuous. In this case, by definition, since $x_n$ strongly converges to $x$, we have that $\nabla f(x_n)$ weakly converges to $\nabla f(x)$. It follows that $||\nabla f(x)||\leq \liminf _{n\rightarrow\infty}||\nabla f(x_n)||$. Hence, if we have $\lim _{n\rightarrow\infty}||\nabla f(x_n)||=0$, it follows also that $\nabla f(x)=0$.)     

The interest in Class $(S)_+$ lies in that one can define topological degree for functions in this class, and  use it to obtain lower bounds on the number of critical points. Moreover, Class $(S)_+$ includes many common PDE, such as Leray-Schauder's where $Id-\nabla f$ is compact. 

\subsection{Preliminaries on measure theory and dynamics on Banach spaces} In the Euclidean setting, a way to characterise "randomness" of a point is to postulate that it is of Lebesgue measure 0. In the Banach space setting, there is a notion which has the same effect, that of "shyness" \cite{christensen}\cite{hunt-etal} which we now introduce. 

{\bf Shyness}. Let $X$ be a Banach space and $S\subset X$ be a Borel set. We say that $S$ is shy if there is a measure $\mu$ on $X$ with the following two requirements: 

i) There is a compact set $A$ such that $0<\mu (A)<\infty$. 

ii) For all $v\in V$, we have $\mu (S+v)=0$. 

The following properties of shyness is useful for the purpose of this paper, see \cite{hunt-etal} for details. 

\begin{proposition} a) Shyness is invariant under translations. 

b) A countable union of shy sets is again a shy set. 

c) If $V$ is a Euclidean space, then a set $S$ is shy iff it has Lebesgue measure $0$.  
\label{PropositionShyness}\end{proposition}

Besides being invariant under translations, sets of Lebesgue measures $0$ on Euclidean spaces also enjoy a stronger invariant property, as stated in Lemma \ref{LemmaMeasureEstimate}. For the purpose of this paper, we also need this stronger invariant property in the infinite dimensional setting. However, it is unclear if the shyness, as defined above, would satisfy this.  Thus, we propose the following slight generalisation of shyness (and also of the local shyness, see Definition 5 and Fact 5 in  \cite{hunt-etal}), which is also invariant under countable unions. 

{\bf Weaker shyness}. Let $X$ be a Banach space and $S\subset X$ be a Borel set.  We say that $S$ is weakly shyness if there are countably many open sets $\{W_i\}_{i\in I}\subset X$ satisfying the following two requirements: 

i) $S=\bigcup _{i\in I}(S\cap W_i)$. 

ii) For each $i\in I$, there are a continuous map $H_i:W_i\rightarrow X$ which has finite distortion (i.e. there is $\lambda _i>0$ so that for all $x,y\in W_i$ we have $||H_i(x)-H_i(y)||\geq \lambda _i||x-y||$), whose image $H_i(W_i)$ is an open subset of $X$, and a shy set $S_i\subset X$ such that $S\cap W_i\subset H_i^{-1}(S_i)$. 

{\bf Remark.} In the above definition, we do not require that $W_i$'s are pairwise distinct. Unlike the case of finite dimension, the requirement that $H_i(W_i)$ is an open subset in ii) above is necessary in general, otherwise we have trivial pathological cases like $V$ has a countable basis $\{e_j\}_{j=1,2,\ldots}$, and $H_i(e_j)=e_{j+1}$ for $j=1,2,\ldots $, and $S_i=$ the subspace generated by $e_2,e_3,\ldots $

We end this subsection with some results on stable-centre manifolds theorem on Banach spaces. We first recall some crucial facts about linear operators on Banach spaces and their spectrum.   

Let $T:X\rightarrow X$ be a bounded linear operator. Then its spectrum $\sigma (T)$ is the set of (complex) numbers $\lambda$ for which $\lambda I-T$ does not have an inverse that is a bounded linear operator. Eigenvalues of $T$ belong to $\sigma (T)$, but $\sigma(T)$ can contain other elements and may be not discrete. However, the spectrum is always a non-empty  compact subset of $\mathbf{C}$. For more detail, please see \cite{kato}.  

In the case of finite dimensions, the Hessian matrix of a real function is symmetric and has only real eigenvalues. In the infinite dimension case, the corresponding notion is self adjointness, and can be defined for Hilbert spaces, or more generally for separable Banach spaces, see e.g. \cite{gill}. 

Notions of saddle points also need more care than the case of finite dimension settings.

{\bf Definition.} Let $X$ be Banach space and $f:X\rightarrow \mathbf{R}$ a $C^2$ map. 

Saddle point: A point $x^*$ is a saddle point of $f$ if it is a critical point (i.e. $\nabla f(x^*)=0$) and $0$ does not belong to the spectrum $\sigma (\nabla ^2f(x^*))$. 

Generalised saddle point: A point $x^*$ is a generalised saddle point of $f$ if it is a critical point ( (i.e. $\nabla f(x^*)=0$)), and  $0$ is not in the closure of $\sigma (\nabla ^2f(x^*))\cap (-\infty, 0)$. 

We have a similar Stable-Center manifold theorem in the Banach space setting, which is crucial for proving local or global avoidance of saddle points for certain iterative algorithms. The following result is also known as Left Stable manifold theorem in the literature. 

\begin{theorem} Let $X=X_+\times X_{-}$ be a Banach space, where $X_{-}$ is of positive dimension. Let $U_+\subset X_+$ and $U_{-}\subset X_{-}$ be open neighbourhoods of $0$, and $f:U_+\times U_{-}\rightarrow X_+\times X_{-}$ be a $C^k$ map ($k\geq 1$), such that $f(0)=0$. Assume that both $V_+=X_+\times \{0\}$ and $V_{-}=\{0\}\times X_{-}$ are invariant by the Jacobian  $Jf(0)$ at $0$. Define $g_+=Jf(0)|_{V_+}:V_+\rightarrow V_+$ and $g_{-}=Jf(0)|_{V_{-}}:V_{-}\rightarrow V_{-}$. Assume moreover that $g_{-}$ is invertible and there are $0<b_1\leq 1<b_2$ such that: $||g_+||\leq b_1$, $b_2\leq ||g_{-}||$. Then there is a small open neighbourhood $U$ of $0$ in $X$ such that the set 
$$W_b=\{x\in U: f^n(x)\in U  \mbox{ for all } n=0,1,2,\ldots \}$$ 
is contained in the graph of a $C^k$ map $h:U\cap X_{+}\rightarrow U\cap X_{-}$. 
\label{TheoremBanachStableCenterManifold}\end{theorem}

In the case $b_1<1<b_2$ and $f$ is a local diffeomorphism near $0$, then this is the Theorem in \cite{irwin}. The proof therein can be adapted for this generalisation. This local Stable-Center manifold turns out to be small, more precisely we have:

\begin{lemma} Let assumptions be as in Theorem \ref{TheoremBanachStableCenterManifold}. Then the graph $\Gamma _h$ of the map $h:U\cap X_{+}\rightarrow U\cap X_{-}$ in the conclusion of the theorem is shy. 
\label{LemmaLocalStableCenterManifoldShy}\end{lemma} 
\begin{proof}
Since $X_{-}$ has positive dimension, we can choose a line $L\subset X_{-}$. Consider $\mu$ the Lebesgue measure on $\{0\}\times L$. Then for all $v=(v_1,v_2)$, the intersection between $\Gamma _h$ and $\{0\}\times L$ is at most 1 point, and hence has measure $0$. Therefore, $\Gamma _h$ is shy, by definition. 
\end{proof}

In the next subsections, we will present the main results concerning Local Backtracking GD, Backtracking GD and (Backtracking) New Q-Newton's method in the Banach space setting. More precisely, we show that Local Backtracking GD solves both Tasks 1 and 2 from the introduction, while BNQN solves Task 1 and locally solves Task 2 (i.e. local Stable-center manifold exists for the dynamics of BNQN near saddle points, and BNQN has quadratic convergence rate near non-degenerate local minima). 

\subsection{Banach Local Backtracking GD and Banach Backtracking GD}

Here we state one concise version of our main results concerning Local Backtracking GD in the Banach space setting.  We note that parts 1, 2 and 3 of Theorem \ref{TheoremMain1} also hold for Backtracking GD. As in the case of finite dimensions, we expect that Backtracking GD also satisfy part 4 of Theorem \ref{TheoremMain1}.

We assume that there is given a canonical isomorphism between $X$ and its dual $X^*$, for example when $X$ is a Hilbert space (see Proposition 8 in \cite{browder}, discussed in comments below, for how to deal with in the case $X$ is not Hilbert).  In the comments afterwards, we will discuss about extensions which are more complicated to state. The definition of Local Backtracking GD in the Banach space setting is straightforwardly extended from its Euclidean space version.

\begin{theorem} Let $X$ be a reflexive Banach space and $f:X\rightarrow \mathbb{R}$ be a $C^2$ function which satisfies Condition C. Moreover, we assume that for every bounded set $S\subset X$, then $\sup _{x\in S}||\nabla ^2f(x)||<\infty$. We choose a point $x_0\in X$ and construct by the Local Backtracking GD  procedure (in the infinite dimensional setting) the sequence $x_{n+1}=x_n-\delta (x_n)\nabla f(x_n)$. Then we have: 

1) Every cluster point of $\{x_n\}$, in the {\bf weak} topology, is a critical point of $f$. 

2) Either $\lim _{n\rightarrow\infty}f(x_n)=-\infty$ or $\lim _{n\rightarrow\infty}||x_{n+1}-x_n||=0$. 

3) Here we work with the weak topology. Let $\mathcal{C}$ be the set of critical points of $f$. Assume that $\mathcal{C}$ has a bounded component $A$. Let $\mathcal{B}$ be the set of cluster points of $\{x_n\}$. If $\mathcal{B}\cap A\not= \emptyset$, then $\mathcal{B}\subset A$ and $\mathcal{B}$ is connected.     

4) Assume that $X$ is separable. Then for generic choices of $\alpha ,\beta ,\delta _0$, there is an exceptional set $\mathcal{E}$ which is weakly shy such that if the initial point $x_0$ is in $X\backslash \mathcal{E}$ and if the sequence $\{x_n\}$ converges - in the {\bf weak} topology, then the limit point cannot be a (generalised) saddle point. 
\label{TheoremMain1}\end{theorem}

 \begin{proof} By the assumption on the Hessian of $f$, we see that $f$ satisfies the conditions needed to apply the Local Backtracking GD procedure, where one define $r(x)=1+||x||$ and $L(x)=\sup _{z\in B(x,r(x))}||\nabla ^2f(z)||$. Let $\{x_n\}$ be a sequence constructed by  the Local Backtracking GD procedure.
 
 1) By the usual estimates in finite dimensional Backtracking GD,  that whenever $K\subset X$ is a bounded set then $\inf _{x\in K}\hat{\delta}(x)>0$, we have that if $\{x_{n_j}\}$ w-converges to $x$, then $\lim _{j\rightarrow\infty}||\nabla f(x_{n_j})||=0$. Then by Condition C, we have that $\nabla f(x)=0$, as wanted. 
 
 2) This follows as in \cite{truong-nguyen}, by using that $\hat{\delta}(x_n)$ is bounded from above. 
  
3) By  Condition C, it follows that $\mathcal{B}\subset \mathcal{C}$ and that $\mathcal{C}$ is w-closed.  We let $X_0\subset X$ be the closure of the vector space generated by $\{x_n\}$. Then $X_0$ is also reflexive, and moreover, it is separable. Moreover, $\mathcal{B}\subset X_0$, and hence $\mathcal{B}\subset \mathcal{C}\cap X_0$. Therefore, from now on, we can substitute $X$ by $X_0$, and substitute $\mathcal{C}$ by $\mathcal{C}\cap X_0$.  We define $\mathcal{C}'=\mathcal{C}\backslash A$. 

We define $\mathcal{B}_0=\mathcal{B}\cap A$, which is a non-empty w-closed set. We will show first that if $z\in \mathcal{B}$ then $z\in \mathcal{B}_0$. We note that $A,\mathcal{B}_0$ are compact in w-topology, and $\mathcal{C}'$ is closed in w-topology, and $A\cap \mathcal{C}'=\emptyset$. Recall that for each $z\in X$, the w-topology has a basis for open neighbourhoods of $z$ the sets of the form $\{x:~||h(x)-h(z)||<\epsilon\}$, where $\epsilon >0$ is a constant and $h:X\rightarrow \mathbb{R}^m$ is a bounded linear map to a finite dimensional space $\mathbb{R}^m$. Therefore, we can find an $\epsilon >0$ and a bounded linear map $h:X\rightarrow \mathbb{R}^m$ so that: for all $x'\in \mathcal{C}'$ and $x\in A$, then $|h(x)-h(z)|>\epsilon$. In particular, $h(A)\cap h(\mathcal{C}')=\emptyset$. Hence, $h(A)$ is a bounded component of $h(\mathcal{C})$. 

By 2) we have that $\lim _{n\rightarrow\infty}||x_{n+1}-x_n||=0$, and hence the same is true for the sequence $\{h(x_n)\}\subset \mathbb{R}^m$. Therefore, by the arguments in \cite{truong-nguyen}, by using the real projective space $\mathbb{P}^m$, we have that the closure in $\mathbb{P}^m$ of the set of  cluster points of $\{h(x_n)\}$  is connected. Note also that the cluster points of $\{h(x_n)\}$ is exactly $h(\mathcal{B})$, is contained in $h(\mathcal{C})$ and has a non-empty intersection with $h(A)$. It follows that $h(\mathcal{B})$ must be contained in $h(A)$ and is connected. In particular, $\mathcal{B}\cap \mathcal{C}'=\emptyset$, as claimed. 

Now to finish the proof, we will show that $\mathcal{B}$ is itself connected.  This can in fact be done as above, by using one connected component of $\mathcal{B}$ in the place of $\mathcal{B}$. 

4) Let $z^*$ be a (generalised) saddle point of $f$. We will show that there is an exceptional set $\mathcal{E}_{z^*}$ such that if the initial point $z_0$ belongs to $X\backslash \mathcal{E}$, then the constructed sequence $\{z_n\}$ cannot converge to $z^*$. 

We first show the existence of a local Stable-Center manifold for the dynamics of Local Backtracking GD in the Banach space setting. As in the case of Riemannian Local Backtracking GD, since the choice of $\alpha$, $\beta$ and $\delta _0$ are random, near $z^*$ the $\widehat{\delta}(x)$ is a constant $\delta$. Therefore, the map $H(x)=x-\widehat{\delta}(x)\nabla f(x)=x-\delta \nabla f(x)$ is $C^1$. Moreover, the map $x\mapsto \widehat{\delta}(x)\nabla f(x)$ is a contraction in the ball $B(z^*,r(z^*))$, and hence both $H(x)$ and $JH(0)$ are invertible. Note that the Jacobian of $H(x)$ is self-adjoint, and the spectrum of $JH(z^*)$ belongs to the set $[0,+\infty )$. Since $z^*$ is a (generalised) saddle point, there are $0<b_1\leq 1<b_2$ such that the spectrum of $JH(z^*)$ belongs to $[0,b_1]\cup [b_2,\infty )$. We let $X_+$ be the subspace of $X$ corresponding to the part $[0,b_1$ of the spectrum, and $X_{-}$ be the subspace corresponding to the part $[b_2,\infty )$ of the spectrum (see Section 4, Chapter III in \cite{kato}). Then we can check that conditions to apply Theorem \ref{TheoremBanachStableCenterManifold} are satisfied, and hence a local Stable-Center manifold $W(z^*)$ exists near $z^*$. 

 By Lemma \ref{LemmaLocalStableCenterManifoldShy}, $W(z^*)$ is shy. Now, as in the Riemannian manifold setting, near each point $y\in X$, the dynamics of Local Backtracking GD in $B(y,r(y))$ is one of a finite collection of maps $x\mapsto H_y(x)=x-\delta \nabla f(x)$, where $\delta$ belongs to a finite set (depending on $y$) and the map $x\mapsto \delta \nabla f(x)$ is a contraction map. Therefore, the map $H_y(x)$ has finite distorsion and is invertible. Hence, the preimage of any weakly shy set is also weakly shyset. 

Since $X$ is reflexive and separable, it follows that $X^*$ is also separable. Then (see point iii) in \cite{dance-sims}) it follows that any bounded subset of $X$, in the weak topology, is metrizable. Also, since $X$ is separable, it follows that $X$ is hereditarily Lindel\"of. That is, any open cover of a subspace of $X$, in the weak topology, has a countable open subcover. We can then proceed as in the proof in the Riemannian case to finish the proof of part 4) here.  
\end{proof}

{\bf Remarks:}

1) For general Banach spaces,  we cannot regard $\nabla f(x_n))\in X^*$ as an element of $X$, and hence the update rule in the theorem is not legit. However,  in parts 1, 2 and 3 of the theorem, we can choose $x_{n+1}=x_n-\hat{\delta}(x_n)v(x_n)$ provided $v(x_n)\in X$ is chosen such that Armijo's condition is satisfied. To this end, we choose (by Hahn-Banach theorem), under the assumption that $X$ is reflexive, a point $v(x_n)\in X$ so that $<\nabla f(x_n),v(x_n)>=||\nabla f(x_n)||^2$ and $||v(x_n)||=||\nabla f(x_n)||$. Hence we can also work with Frechet differentiation as in \cite{blank-rupprecht}. Likewise, we can use other versions of Backtracking GD, under less restrictions on the function $f$ ($f$ being $C^1$ is enough).

Note that as of current, for parts 4 and 5 of the theorem we need to use Stable-Central manifold, and hence if the map $x\mapsto x-\hat{\delta}(x)v(x)$ is not $C^1$ near saddle points, then it is difficult to proceed, as of current. See however point 5 below for an idea on using duality mappings to deal with this in general. 

2) This theorem can be extended to several other versions of Backtracking GD, including Unbounded Backtracking GD. The conclusion of parts 4 and 5 of the theorem are also valid if we use the strong topology instead of weak topology.

3) In part 3) of Theorem \ref{TheoremMain1}, if $f$ has at most countably many critical points (including Morse functions), then the conclusion is simple to state: Either $\lim _{n\rightarrow\infty}||x_n||=\infty$ or $\{x_n\}$ w-converges to a critical point of $f$. 

4) Duality mapping: 

If we are willing to redefine the norms on $X$, then $v(x_n)$ (in point 1 above) is unique and depends continuously on $x_n$. The argument presented here follows Proposition 8 in \cite{browder}. Assume that $X$ is reflexive. Then we can renorm both $X$ and $X^*$ so that both of them are locally uniformly convex. Then, we have a well-defined map $J:X^*\rightarrow X$ with the following properties: $<y,J(y)>=||y||^2$ and $||J(y)||=||y||$ for all $y\in X^*$. This map is called normalized duality mapping, and it is bicontinuous. Moreover, it is monotone, that is $<y_1-y_2,J(y_1)-J(y_2)>\geq 0$, for all $y_1,y_2\in X^*$. Indeed, 
\begin{eqnarray*}
<y_1-y_2,J(y_1)-J(y_2)>&=& [(||y_1||-||y_2||)^2] + [||y_1||\times ||J(y_2)||-<y_1,J(y_2)>] \\
&+& [||y_2||\times ||J(y_1)||-<y_2,J(y_1)>],
\end{eqnarray*}
and each square bracket is non-negative. Hence, the map $J$ is in Class $(S)_+$. 

Our GD update rule is now: $x\mapsto x-\delta (x)J(\nabla f(x))$, which is continuous near critical points. Even though it is not $C^1$ in general, it seems good enough that conclusions of  parts 4 and 5 of Theorem \ref{TheoremMain1} should follow. For example, if $X^*$ is the space $l^p$ ($1<p<\infty$), i.e. the space of sequences $y=(y^{(1)},y^{(2)},\ldots )$ so that $\sum _{j}|y^{(j)}|^p<\infty$, with norm $||y||=(\sum _j|y^{(j)}|^p)^{1/p}$, then the normalized duality mapping is: 
\begin{eqnarray*}
J(y)=\frac{1}{||y||^{p-2}}(|y^{(1)}|^{p-1}sign (y^{(1)}),|y^{(2)}|^{p-1}sign (y^{(2)}),\ldots )
\end{eqnarray*}
Assume that $p\geq 2$. The map $J(y)$ is $C^1$ outside of the point $y=0$, and its derivative is locally bounded near this point. (To see this, the readers can simply work out the differentiability of the following map in $2$ variables $|s|^{p-2}s/(|s|^p+|t|^p)^{(p-2)/p}$.)  Therefore, in this case, the GD dynamical systems $x\mapsto \mathcal{F}(x)=x-\delta (x)J(\nabla f(x))$ is $C^1$ in small punctured neighbourhoods of critical points of $f$.   

\subsection{Banach (Backtracking) New Q-Newton's method} In this subsection, we work under assumptions on the Banach space $X$ as in the previous subsection. That is, we assume that there is given a canonical isomorphism between $X$ and its dual $X^*$, for example when $X$ is a Hilbert space (and see the discussion in the end of the previous subsection for how to deal with in the case $X$ is not Hilbert).

 In the definition of (Backtracking) New Q-Newton's method in the finite dimensional setting, there are two While loops, where the first While loop is to choose $\delta _j$ for which all the eigenvalues of the matrix $\nabla^2f(x_n)+\delta_j \|\nabla f(x_n)\|^{1+\alpha}Id$ have absolute value $\geq $ $\|\nabla f(x_n)\|^{1+\alpha}\inf _{i\not=j } |\delta _i-\delta _j|/2$. 

A straightforward of this definition to the infinite dimensional setting is as follows. We choose a discrete sequence $\{\delta _k\}$. We then require that all elements of the spectrum of $\nabla^2f(x_n)+\delta_j \|\nabla f(x_n)\|^{1+\alpha}Id$ have absolute value $\geq $ $\|\nabla f(x_n)\|^{1+\alpha}\inf _{i\not=j } |\delta _i-\delta _j|/2$. Then obviously, the operator $\nabla^2f(x_n)+\delta_j \|\nabla f(x_n)\|^{1+\alpha}Id$ is invertible (under the  assumption that $\nabla f(x_n)\not= 0$). 

(The second While loop, in the definition of BNQN, is kept the same.)

While the above extension works theoretically, it poses practical challenges in the infinite dimensional setting, and we propose the following version which is simpler and does not require performing the first While loop. In the finite dimensional setting, it is simpler than the previous version. First, let $T$ be a self adjoint bounded linear operator over $X$, and let $V_+\subset X$ be the subspace corresponding to the part of $\sigma (T)\cap [0,+\infty )$ and $V_{-}$ be the subspace corresponding to the part of   $\sigma (T)\cap (-\infty , 0]$.  We let $pr_{T,+}$ and $pr_{T,-}$ be the natural projections from $X$ to $V_+$ and $V_{-}$ (see Section 4, Chapter 3 in \cite{kato}).  We then define $|T|$ to be the linear operator given by: 
$$Abs(T)(x)=pr_{T,+}(Tx)-pr_{T,-}(Tx).$$

It follows that for all $\delta >0$, if $\nabla f(x_n)\not= 0 $, then $Abs(\nabla^2f(x_n))+\delta \|\nabla f(x_n)\|^{1+\alpha}Id$ is invertible.

\medskip
{\color{blue}
 \begin{algorithm}[H]
\SetAlgoLined
\KwResult{Find a minimum of $F:\mathbb{R}^k\rightarrow \mathbb{R}$}
Given: $\delta >0$\  and $\alpha >0$;\\
Initialization: $x_0\in \mathbb{R}^k$\;
 \For{$n=0,1,2\ldots$}{

$A_n:=Abs(\nabla^2f(x_n))+\delta \|\nabla f(x_n)\|^{1+\alpha}Id$\\
$v_n:=A_n^{-1}\nabla f(x_n)=pr_{A_n,+}(v_n)+pr_{A_n,-}(v_n)$\\
$w_n:=pr_{A_n,+}(v_n)-pr_{A_n,-}(v_n)$\\
When $f$ does not have compact sublevels, normalise $w_n:=w_n/\max\{1,||w_n||\}$\\

$\gamma :=1$\\

 \While{$f(x_n-\gamma w_n)-f(x_n)+\frac{\gamma}{3}\nabla f(x_n).w_n>0 $}{$\gamma =\gamma /2$}

$x_{n+1}:=x_n-\gamma w_n$
   }
  \caption{Simplified Backtracking New Q-Newton's method  (Simplified BNQN)} \label{table:alg}
\end{algorithm}
}
\medskip

The proofs for the properties of  the Riemannian Backtracking New Q-Newton's method and Banach Local Backtracking GD can be combined to attain the following result. 

\begin{theorem} Let $X$ be a reflexive Banach space and $f:X\rightarrow \mathbb{R}$ be a $C^2$ function which satisfies Condition C. Moreover, we assume that for every bounded set $S\subset X$, then $\sup _{x\in S}||\nabla ^2f(x)||<\infty$. We choose a point $x_0\in X$ and construct by Simplified BNQN. Then we have: 

1) Every cluster point of $\{x_n\}$, in the {\bf weak} topology, is a critical point of $f$. 

2) Either $\lim _{n\rightarrow\infty}f(x_n)=-\infty$ or $\lim _{n\rightarrow\infty}||x_{n+1}-x_n||=0$. 

3) Here we work with the weak topology. Let $\mathcal{C}$ be the set of critical points of $f$. Assume that $\mathcal{C}$ has a bounded component $A$. Let $\mathcal{B}$ be the set of cluster points of $\{x_n\}$. If $\mathcal{B}\cap A\not= \emptyset$, then $\mathcal{B}\subset A$ and $\mathcal{B}$ is connected.     

4) Let $z^*$ be a saddle point of $f$. Then there exists a local Stable-center manifold $W(z^*)$ (which is shy) for the dynamics of Simplified BNQN near $z^*$. 

5) If $z^*$ is a non-degenerate local minimum of $f$, then Simplified BNQN has local quadratic rate of convergence near $z^*$. 
\label{TheoremMain2}\end{theorem}

We note that to prove the (full) global avoidance of saddle points for Simplified BNQN in the Banach space setting, we need to check that locally the dynamics of Simplified BNQN is one of a finite collection of invertible maps with finite distortion. This imposes serious challenges, given that the spectrum of a linear operator in infinite dimension is rather complicated, and maybe not discrete. To be able to show global avoidance of saddle points for BNQN in the Banach space setting, it is necessary to extend Theorem \ref{TheoremRandomnessLambda}  to the infinite dimensional setting, which is beyond the current paper.

\end{document}